\documentclass[11pt,notitlepage,reqno]{amsart}
\usepackage[margin=1in]{geometry}
\usepackage{mathtools} 
\usepackage{amssymb} 
\usepackage{mathrsfs}
\usepackage{bbm}
\usepackage{enumitem}
\usepackage{setspace}
\usepackage[dvipsnames]{xcolor}
\usepackage{multicol}
\usepackage{lastpage}
\usepackage{float}
\usepackage{graphicx}
\usepackage{subfigure}
\usepackage{outlines}
\usepackage[labelfont=small,sc]{caption}
\usepackage{hyperref}
\usepackage[noabbrev,capitalise]{cleveref}
\usepackage{bookmark}
\bookmarksetup{
	numbered,
	open
}

\usepackage{amsthm, thmtools}
\newtheorem{theorem}{Theorem}[section]
\newtheorem*{theorem*}{Theorem}
\newtheorem{lemma}[theorem]{Lemma}
\newtheorem{proposition}[theorem]{Proposition}
\newtheorem*{proposition*}{Proposition}
\newtheorem{corollary}[theorem]{Corollary}
\newtheorem*{corollary*}{Corollary}

\newtheorem{definition}[theorem]{Definition}

\newtheorem*{problem*}{Problem}

\theoremstyle{remark}
\newtheorem{remark}[theorem]{Remark}

\newtheorem{example}[theorem]{Example}

\newtheorem*{claim*}{Claim}

\numberwithin{equation}{section}
\numberwithin{theorem}{section}
\usepackage{chngcntr}

\allowdisplaybreaks
\usepackage{tikz}
\usetikzlibrary{matrix,chains,shapes,arrows,scopes,positioning,patterns}
\usetikzlibrary{decorations.markings,decorations.pathreplacing}

\newcommand{\Irr}{\operatorname{Irr}}
\newcommand{\Lin}{\operatorname{Lin}}
\newcommand{\sgn}{\operatorname{sgn}}
\newcommand{\Syl}{\operatorname{Syl}}
\newcommand{\down}{\big\downarrow}
\newcommand{\triv}{\mathbbm{1}}
\newcommand{\C}{\mathbb{C}}
\newcommand{\N}{\mathbb{N}}

\newcommand{\Z}{\mathbb{Z}}
\newcommand{\cA}{\mathcal{A}}
\newcommand{\cB}{\mathcal{B}}

\newcommand{\cD}{\mathcal{D}}
\newcommand{\cL}{\mathcal{L}}
\newcommand{\cP}{\mathcal{P}}
\newcommand{\sh}{\mathsf{h}}

\newcommand{\sL}{\mathsf{L}}
\newcommand{\sX}{\mathsf{X}}
\newcommand{\sY}{\mathsf{Y}}

\newcommand{\blue}{\color{black}}

\newcommand{\sbc}{\mathcal{Z}}
\newcommand{\set}[1]{\{#1\}}
\newcommand{\Hook}{\mathcal{H}}
\newcommand{\Ahook}{\mathcal{AH}}
\newcommand{\ahook}[2]{\mathsf{ah}_{#1}(#2)}
\newcommand{\hook}[2]{\mathsf{h}_{#1}(#2)}
\newcommand{\Chi}{\mathcal{X}}
\newcommand{\ind}{\big\uparrow}
\newcommand{\res}{\big\downarrow}
\newcommand{\inprod}[1]{\left\langle#1\right\rangle}
\newcommand{\abs}[1]{\big\lvert#1\big\rvert}
\newcommand{\inv}{^{-1}}

\newcommand{\Part}{\cP}
\newcommand{\lin}{\cL}

\begin{document}

\title[Minimal numbers of linear constituents in Sylow restrictions for $S_n$]{Minimal numbers of linear constituents in Sylow restrictions for symmetric groups}

\author{Bim Gustavsson}
\address[B.~Gustavsson]{School of Mathematics, Watson Building, University of Birmingham, Edgbaston, Birmingham B15 2TT, UK}
\email{bxg281@bham.ac.uk}

\author{Stacey Law}
\address[S.~Law]{School of Mathematics, Watson Building, University of Birmingham, Edgbaston, Birmingham B15 2TT, UK}
\email{s.law@bham.ac.uk}

\begin{abstract}
	Let $p$ be any prime. We determine precisely those irreducible characters of symmetric groups which contain at most $p$ distinct linear constituents in their restriction to a Sylow $p$-subgroup, answering a question of Giannelli and Navarro in \cite{GN}. Moreover, we identify all of the linear constituents of such characters, and in the case $p=2$ explicitly calculate a new class of Sylow branching coefficients for symmetric groups indexed by `almost hook' partitions.
\end{abstract}

\maketitle

\section{Introduction}

A major theme of research in the representation theory of finite groups is to investigate the relationship between the characters of a finite group $G$ and those of its Sylow subgroups. In particular, the restriction of characters of $G$ to Sylow $p$-subgroups of $G$ and counting the number of linear constituents in these has been studied in connection with the McKay Conjecture (recently a theorem), for instance in \cite{NTV}, \cite{G} and \cite{GN}. In this article, we continue this study into linear constituents of character restrictions to Sylow subgroups in the case of symmetric groups.

Let $n\in\N$. It is well known that the set $\Irr(S_n)$ of irreducible characters of the symmetric group $S_n$ is in natural bijection with $\cP(n)$, the set of partitions of $n$. Given $\lambda\in\cP(n)$, we denote by $\chi^\lambda$ the corresponding character in $\Irr(S_n)$. 
For any prime $p$ and $P\in\Syl_p(S_n)$, it was shown in \cite[Theorem A]{GN} that if $\chi\in\Irr(S_n)$ and $p\mid\chi(1)$, then the restriction $\chi\down_P$ contains at least $p$ different linear constituents.

Our first main result is an explicit description of those $\chi\in\Irr(S_n)$ which attain this lower bound, answering a question posed in \cite{GN}. In fact, we determine all those $\chi$ with \textit{at most} $p$ distinct linear constituents, thereby determining those $\chi$, irrespective of whether their degree is divisible by $p$, which have minimal numbers of linear constituents upon restriction to a Sylow subgroup of $S_n$. 
Given a character $\phi$ of a finite group $G$, we let $\Lin(\phi)$ denote the set of linear constituents of $\phi$.
We separate our results for $p=2$ from those for odd primes $p$, since in order to state the classification in the latter case it will be useful to introduce some additional notation.

\begin{theorem}\label{thm:main-p=2}
	Let $n\in\N$ and $\lambda\in\cP(n)$. Let $P\in\Syl_2(S_n)$. Then $|\Lin(\chi^\lambda\down_P)|=2$ if and only if one of the following holds:
	\begin{enumerate}
		\item $n=2^k$ for some $k\in\N_{\ge 2}$ and $\lambda=(2^{k-1},2,1^{2^{k-1}-2})$.
		\item $n=2^k+1$ for some $k\in\N$ and $\lambda=(n-t,1^t)$ for some $1\le t\le n-2$,
		\item $n=8$ and $\lambda\in\{(5,3),(3^2,2),(2^3,1^2)\}$, or $n=9$ and $\lambda=(3^3)$.
	\end{enumerate}
	Moreover, in each of the above cases, $\inprod{\chi^\lambda\down_P,\psi}=1$ for all $\psi\in\Lin(\chi^\lambda\down_P)$.
\end{theorem}

We remark that when $p=2$ it is already known when 
$|\Lin(\chi^\lambda\down_P)|=1$ holds from \cite[Theorem 1.2]{G}; this is recorded in \cref{sec:ah} below.

For any $n,t\in\N$, we let $\cB_n(t)$ denote the set of partitions of $n$ whose first part and number of parts are both at most $t$. Below we consider only $n\ge p$ since otherwise $P\in\Syl_p(S_n)$ is trivial, and only $\lambda\notin\{(n),(1^n)\}$ as otherwise $\chi^\lambda\res_P$ is trivial.

\begin{theorem}\label{thm:main-odd-p}
	Let $p$ be an odd prime and suppose $n\ge p$ is a natural number. Let $\lambda\in\cP(n)\setminus\set{(n),(1^n)}$ and $P\in\Syl_p(S_n)$.
	\begin{enumerate}
		\item Suppose $n=p^k$ for some $k\in\N$. Then
		\[ |\Lin(\chi^\lambda\down_P)| = \begin{cases}
			p-1 & \text{if } \lambda\in\{(n-1,1),(2,1^{n-2})\},\\
			p & \text{if }\lambda\in\cB_n(n-2)\setminus\cB_n(n-p),\text{ or }(p,k)=(3,2)\text{ and }\lambda=(3^3),
		\end{cases} \]
		and $|\Lin(\chi^\lambda\down_P)|>p$ otherwise.
		
		\item Suppose $n$ is not a power of $p$. If $n=p^k+i$ for some $k\in\N$ and $i\in\set{1,2,\dots,p-1}$, then 
			\[ \abs{\Lin(\chi^\lambda\res_P)} = \begin{cases}
				{\blue p-1} & {\blue \text{if } (p,n,\lambda)=(3,4,(2^2)),}\\
				p &\text{if } \lambda\in\cB_n(n-1)\setminus\cB_n(n-p) {\blue \text{ and } (p,n,\lambda)\ne(3,4,(2^2)),}\\
				{\blue p} & {\blue \text{if } k=1 \text{ and }\lambda\in\cB_n(n-p).} 
			\end{cases} \]
		In all other cases, $\abs{\Lin(\chi^\lambda\res_P)}>p$.
	\end{enumerate}
\end{theorem}

For all primes $p$, we identify all of the distinct linear constituents when $|\Lin(\chi^\lambda\down_P)|\le p$ in \cref{rem:identify-two-linears,rem:identify-few-linears-odd} below, after a parametrisation of the relevant linear characters of Sylow subgroups of symmetric groups has been fixed.  

Theorems~\ref{thm:main-p=2} and~\ref{thm:main-odd-p} can also be thought of as results on the positivity of the so-called Sylow branching coefficients corresponding to linear characters of Sylow subgroups of symmetric groups. Given a finite group $G$, a prime number $p$ and $P\in\Syl_p(G)$, for each $\chi\in\Irr(G)$ we can decompose
\[ \chi\down_P = \sum_{\phi\in\Irr(P)} \sbc^\chi_\phi\phi, \]
where the quantities $\sbc^\chi_\phi$ are termed \textit{Sylow branching coefficients}. Such coefficients describe character restriction and induction between finite groups and their Sylow subgroups. When $G=S_n$ and $\chi=\chi^\lambda$ for some $\lambda\in\cP(n)$, we often abbreviate $\sbc^{\chi^\lambda}_{\phi}$ to $\sbc^\lambda_\phi$.

Part of the classification in \cref{thm:main-p=2} follows from \cref{thm:almost-hook-sbcs} below, which gives us new and precise information on the Sylow branching coefficients for symmetric groups when $p=2$. 
Despite some recent advances such as \cite{LO} and \cite{GV-vietnam}, relatively little is known about such Sylow branching coefficients in the case $p=2$. For instance, although we have asymptotic results such as \cite[Theorem C]{LO}, determining when $\sbc^\lambda_\triv$ is positive where $\triv$ denotes the trivial character of a Sylow 2-subgroup of $S_n$ remains an open problem in general. On the other hand, the positivity of $\sbc^\lambda_\phi$ has been determined for almost all partitions $\lambda$ for all \textit{linear} characters $\phi$ of a Sylow $p$-subgroup of $S_n$ when $p$ is odd \cite{GL2,L}, and recently also all non-linear irreducible $\phi$ when $p\ge 5$ \cite{GL3}.

Focusing instead on a special class of partitions $\lambda$, we are able to describe not just the positivity of $\sbc^\lambda_\phi$ for all linear characters $\phi$ of $P\in\Syl_2(S_{2^k})$, but in fact give their values explicitly in the second main result of this paper. This also generalises \cite[Lemma 3.6]{LO} from the trivial character to all linear characters of $P$.

\begin{theorem}\label{thm:almost-hook-sbcs}
	Let $k\in\N_{\ge 2}$ and $P\in\Syl_2(S_{2^k})$. Suppose $x\in\{0,1,\dotsc,2^k-4\}$ and let the set of linear characters of $P$ be parametrised as $\{\cL(\sh_{2^k}(y)) \mid y\in\{0,1,\ldots,2^k-1 \}\}$ as described in \cref{thm: hook lins} below. Then
	\[ \sbc^{(2^k-2-x,2,1^x)}_{\cL(\sh_{2^k}(y))} = \begin{cases}
		\binom{k-1}{x-B(y)}-1 & \text{if }x\in\{y-1,y-2\},\\
		\binom{k-1}{x-B(y)} & \text{if }x\in \{ B(y), B(y)+1,\ldots, B(y)+k-1\}\setminus\{y-1,y-2\},\\
		0 & \text{otherwise,}
	\end{cases} \]	
	where $B(y)=y-b(y)-1$ and $b(y)$ denotes the sum of the digits in the binary expansion of $y$.
\end{theorem}
 
This article is structured as follows. In \cref{sec:prelims} we recall the relevant background on the representation theory of symmetric groups and their Sylow subgroups. 
In \cref{sec:ah} we compute explicitly the Sylow branching coefficients for symmetric groups $S_{2^k}$ corresponding to $\chi^\lambda$ and $\phi$ for all `almost hook' partitions $\lambda$ and all linear characters $\phi$ of $P\in\Syl_2(S_{2^k})$, proving \cref{thm:almost-hook-sbcs}. It turns out that all other partitions $\lambda\in\cP(2^k)$ (excluding hooks) have $|\Lin(\chi^\lambda\down_P)|>2$; this is shown in the first part of \cref{sec:p=2}. The structure of Sylow $p$-subgroups of $S_n$ depends on the $p$-adic expansion of $n$, and so from our results for $n=2^k$ we are then able to complete the proof of \cref{thm:main-p=2} for arbitrary $n\in\N$ in the second part of \cref{sec:p=2}. In \cref{sec:odd} we address the case of odd primes $p$ and prove \cref{thm:main-odd-p}.

\begin{remark}
	Another way of rephrasing \cref{thm:main-p=2,thm:main-odd-p} is the classification of the set
	\begin{equation}\label{eqn:G}
		\Phi(G,p):=\{\phi\in\Irr(G)\mid |\Lin(\phi\down_P)|=p \}
	\end{equation}
	where $P\in\Syl_p(G)$, for $G=S_n$ and all primes $p$.
	Since the representation theory of $A_n$ is closely related to that of $S_n$, our main results form a first step towards classifying $\Phi(A_n,p)$.
	When $p$ is odd, $\Syl_p(S_n)=\Syl_p(A_n)$, but this is not the case when $p=2$. To determine $\Phi(A_n,2)$ would require further study of Sylow branching coefficients $\sbc^\chi_\phi$ for $S_n$ and $A_n$, including for certain non-linear $\phi\in\Irr(P)$ where $P\in\Syl_2(S_n)$, which will be the topic of future investigation. Initial computations 
	suggest that $\Phi(A_8,2)=\{ \chi^{(5,3)}\down_{A_n} = \chi^{(2^3,1^2)}\down_{A_n}  \}$ and
	\[ \Phi(A_{2^k+1},2)= \big\{ \chi^\lambda\res_{A_n} \mid \lambda \in \{ \hook{2^k +1}{x} \mid x\in \{0,1,...,2^k\} \setminus \{0, 2^{k-1}-1, 2^{k-1}, 2^{k-1}+1, 2^k \} \} \big\} \]
	whenever $k\ge 3$, and $\Phi(A_n,2)=\varnothing$ otherwise. Moreover, when $\Phi(A_n,2)\ne\varnothing$, then $\sbc^\chi_\psi=1$ for all $\chi\in\Phi(A_n,2)$ and $\psi\in\Lin(\chi\down_P)$.
	
	Another interesting question would be to determine $\Phi(G,p)$ and their associated Sylow branching coefficients for other simple groups $G$. For $p=2$ and $P\in\Syl_2(G)$:
	\begin{itemize}
		\item When $G=\operatorname{PSL}(3,3)$, we find that $\Phi(G,2)=\{\chi_1,\chi_2\}$ for some $\chi_1\ne\chi_2\in\Irr(G)$. 
		We note that $\sbc^{\chi_1}_\psi=2,2$ and $\sbc^{\chi_2}_\psi=2,1$ as $\psi$ ranges over $\Lin(\chi_i\down_P)$, unlike in \cref{thm:main-p=2}.
		\item However, $G=\operatorname{PSU}(3,3)$ has a unique $\chi\in\Irr(G)$ s.t.~$\abs{\Lin(\chi\res_P)}=2$, and $\sbc^\chi_\psi=1$ for all $\psi\in\Lin(\chi\down_P)$, like in \cref{thm:main-p=2}.
	\end{itemize}
\end{remark}

\section{Preliminaries}\label{sec:prelims}

Let us begin by fixing some notation which will be used throughout this article. 
Let $\N=\set{1,2,\dots}$ and $\N_0 = \N \cup \set{0}$. For $k\in\N_0$ we let $\N_{\geq k} = \set{k,k+1,k+2,\dots}$. 

Let $n\in\N$, then $\Part(n)$ denotes the set of all partitions of size $n$. Given a partition $\lambda=(\lambda_1, \dots, \lambda_t)$ of size $n$ (meaning $\lambda_1+\cdots+\lambda_t=n$ and written $|\lambda|=n$), we call $\lambda_1,\dots,\lambda_t$ the \textit{parts} of $\lambda$. The number of parts of $\lambda$ is $t$, which we also denote by $\ell(\lambda)$. 

We denote by $\lambda'$ the \textit{conjugate} partition of $\lambda$ and we say that $\lambda$ is \emph{self-conjugate} if $\lambda=\lambda'$. For $\lambda\ne\lambda'$, if $i=\min_j\set{\lambda_j \neq \lambda'_j}$ is such that $\lambda_i > \lambda'_i$, then we say that $\lambda$ is \emph{wide}. If $\lambda'$ is wide, then we say that $\lambda$ is \emph{tall}. Hence any partition is exactly one of wide, tall or self-conjugate. 

The following partition shapes are key concepts that will appear throughout this article.
\begin{definition}
	Let $n\in\N$ and $\lambda\in\Part(n)$.
	\begin{enumerate}[label=\alph*)]
		\item If $\lambda=(n-x,1^x)$ for some $x\in\{0,1,\dotsc,n-1\}$, then we call $\lambda$ a \emph{hook partition}. We set $\hook{n}{x}:=(n-x,1^x)$ and denote by $\Hook(n)$ the subset of $\Part(n)$ consisting of hook partitions.
		\item If $\lambda=(n-2-x,2,1^x)$ for some $x\in\{0,1,\dotsc,n-4\}$, then we call $\lambda$ an \emph{almost hook partition}. We set $\ahook{n}{x}:=(n-2-x,2,1^x)$ and denote by $\Ahook(n)$ the subset of $\Part(n)$ consisting of almost hook partitions.
	\end{enumerate}
\end{definition}

For a finite group $G$, we denote by $\Irr(G)$ the complete set of irreducible complex characters of $G$. Let $\Lin(G)$ denote the subset of $\Irr(G)$ consisting of all linear characters. For a character $\psi$ of $G$, we let $\Lin(\psi) = \set{\phi\in\Lin(G) \mid \inprod{\psi, \phi} > 0}$. For a prime $p$, let $\Irr_{p'}(G)$ denote the set of irreducible characters of $G$ whose degree is coprime to $p$. 

\subsection{The representation theory of the symmetric group}\label{sec:rep-sn}
	We refer the reader to \cite{JK, Sagan} for further detail.
	For $n\in\N$, let $S_n$ denote the symmetric group on $n$ objects. It is well known that $\Irr(S_n)$ is in natural bijection with $\cP(n)$. For $\lambda\in\cP(n)$, we denote by $\chi^\lambda$ the corresponding irreducible character of $S_n$. It is therefore natural to use the same terminology defined for partitions for irreducible characters of $S_n$, e.g.~we say that $\chi\in\Irr(S_n)$ is a wide hook character if $\chi=\chi^\lambda$ for some $\lambda\in\Part(n)$ and $\lambda$ is a wide hook partition. Hence we can, without causing confusion, simply write hook instead of hook character or hook partition. 
	
	It is well known (see \cite[2.1.8]{JK}) that 
	\begin{equation}\label{eqn:sgn}
		\sgn_{S_n}\cdot \chi^\lambda = \chi^{\lambda'},
	\end{equation}
	where $\sgn_{S_n}$ denotes the sign character of $S_n$. Hence $\abs{\Lin(\chi^\lambda\res_{P_n})} = \abs{\Lin(\chi^{\lambda'}\res_{P_n})}$, and so when counting linear constituents of $\chi^\lambda\res_{P_n}$ it is enough to only consider wide or self-conjugate characters $\chi^\lambda$. 
	
	Let $m,n\in\N$ and $\mu\in\Part(m)$, $\nu\in\Part(n)$. For $\lambda\in\Part(m+n)$, the Littlewood--Richardson rule describes how to decompose $\chi^\lambda\res_{S_m \times S_n}$ as a sum of irreducible characters. 
	The following statement of the rule uses the notation of \cite[Theorem 4.9.4]{Sagan}.
	\begin{theorem}\label{thm:LR}
		Let $m,n\in\N$ and let $\lambda\in\Part(m+n)$. Then 
			$$ \chi^\lambda\res_{S_m \times S_n} = \sum_{\mu\in\Part(m),\ \nu\in\Part(n)} c^\lambda_{\mu,\nu}\cdot\chi^\mu\times\chi^\nu,$$
		where $c^\lambda_{\mu,\nu}$ equals the number of semistandard tableaux $T$ such that $T$ has shape $[\lambda\setminus\mu]$ and content $\nu$, and the reverse row word of $T$ is a ballot sequence.
	\end{theorem}
	The coefficients $c^\lambda_{\mu,\nu}$ are called the \emph{Littlewood--Richardson coefficients}. From \eqref{eqn:sgn} and \cref{thm:LR} it is clear that $c^\lambda_{\mu,\nu} = c^\lambda_{\nu,\mu} = c^{\lambda'}_{\mu',\nu'}$. We observe that if $m=n$ and $\lambda\in\Ahook(2n)$, then $c^\lambda_{\mu,\nu}>0$ implies that $\mu\in\Hook(n)$ and $\nu\in\Ahook(n)\cup\Hook(n)$, or $\mu\in\Ahook(n)\cup\Hook(n)$ and $\nu\in\Hook(n)$.

	We note that restrictions to Sylow subgroups of irreducible characters of symmetric groups always have at least one linear constituent.
	\begin{proposition}\label{prop: at least one lin}
		Let $n\in\N$, $p$ be a prime and $P\in\Syl_p(S_n)$. Then $\abs{\Lin(\chi\res_{P})}>0$ for all $\chi\in\Irr(S_n)$. 
	\end{proposition}
	\begin{proof}
		If $p\mid\chi(1)$, then $\abs{\Lin(\chi\res_{P})}\geq p > 1$ by \cite[Theorem A]{GN}. Now suppose that $p\nmid \chi(1)$. Then $\inprod{\chi\res_{P},\psi}>0$ for some $\psi\in\Irr_{p'}(P)$. The statement now follows from the fact that $\Irr_{p'}(P)=\Lin(P)$, since $\phi(1)\mid|G|$ whenever $\phi\in\Irr(G)$.
	\end{proof}

\subsection{Wreath products}

	The following will be important when we want to understand the relation between $S_n$ and its Sylow $p$-subgroups; see \cite[\textsection 4]{JK} for further detail. 
	
	Fix $n\in\N$. Let $G$ be a finite group and let $H\leqslant S_n$. The \emph{wreath product} of $G$ by $H$ is the group $G\wr H$ on the set $\set{(g_1,\dots,g_n;h)\mid g_1,\dots,g_n\in G,\ h\in H}$ with group multiplication $(g_1,\dots,g_n; h)(g'_1,\dots,g'_n; h') = (g_1g'_{h\inv(1)}, \dots, g_ng'_{h\inv(n)}; hh')$. 
	
	Fix a prime $p$ and let $n\in \N$. We denote by $P_n$ a Sylow $p$-subgroup of $S_n$. In particular, $P_1$ is the trivial group and $P_p$ is cyclic of order $p$. For $k\in\N$, we have $P_{p^k} = P_p^{\wr k} = P_p \wr P_p \wr \dots \wr P_p$, i.e.~the $k$-fold wreath product of copies of $P_p$. Letting $n=a_0p^0+a_1p^1+\dots+a_tp^t$ be the $p$-adic expansion of $n\in\N$ (i.e.~$a_i\in\{0,1,\dotsc,p-1\}$ for all $i$), then 
	\begin{equation}\label{eqn:Pn}
		P_n \cong P_{p^0}^{\times a_0}\times P_{p^1}^{\times a_1}\times \dots \times P_{p^t}^{\times a_t}.
	\end{equation}
	
	Next, we need to consider the irreducible characters of wreath products.
	Let $V$ be a left $\C G$--module, affording the character $\psi$. Let $V^{\otimes n}$ denote the $n$-fold tensor product of $V$ over $\C$. The action of $(g_1,\dots,g_n;h)$ on the basis elements of $V^{\otimes n}$ defined by $(g_1,\dots,g_n;h)v_{i_1}\otimes\dots\otimes v_{i_n} = g_1v_{h\inv(i_1)}\otimes \dots \otimes g_nv_{h\inv(i_n)}$, extended linearly, turns $V^{\otimes n}$ into a $\C (G\wr H)$--module. Let $\widetilde{V^{\otimes n}}$ denote this module, and $\widetilde{\psi}$ the character afforded by $\widetilde{V^{\otimes n}}$. Let $\phi$ be a character of $H$, then we also denote by $\phi$ its inflation from $H$ to $G\wr H$, noting that $G^{\times n}\cong\{(g_1,\dotsc,g_n;1)\mid g_i\in G\}\trianglelefteq G\wr H$. That is, $\phi((g_1,\dots,g_n;h)):=\phi(h)$. Hence $\widetilde{\psi}\cdot\phi$ is a character of $G\wr H$, which we denote by $\Chi(\psi;\phi)$.
	Indeed, if $\psi\in\Irr(G)$ then $\Irr(G\wr H\mid\psi^{\times n})$, the set of irreducible characters of $G\wr H$ whose restriction to $G^{\times n}$ contains $\psi^{\times n}$ as a constituent, is precisely
	\[ \{ \Chi(\psi;\phi)\mid \phi\in\Irr(H) \} \]
	(see \cite[Corollary 6.17]{Isaacs}). To ease notation, when the meaning is clear we will sometimes abbreviate $\chi^\lambda$ to simply $\lambda$, such as $\Chi(\mu;\nu):=\Chi(\chi^\mu;\chi^\nu)$.
	
	We record $\Irr(G\wr H)$ when $H=C_p$ for some prime $p$. Each $\chi\in\Irr(G\wr C_p)$ takes one of the following forms:
		\begin{enumerate}[label=\roman*)]
			\item $\chi = \psi_{i_1}\times \dots\times \psi_{i_p}\ind_{G^{\times p}}^{G\wr C_p}$, where $\psi_{i_1},\dots,\psi_{i_p}\in\Irr(G)$ are not all equal; or 
			\item $\chi = \Chi(\psi;\phi)$ for some $\psi\in\Irr(G)$ and $\phi\in\Irr(C_p)$.
		\end{enumerate}
	When i) holds, we have $\chi\res_{G^{\times n}} = \sum_{\sigma\in C_p} \psi_{i_{\sigma(1)}}\times\dots\times\psi_{i_{\sigma(p)}}$, and when ii) holds, we have that $\chi\res_{G^{\times p}} = \psi^{\times p}$. Unless stated otherwise, we let $\Irr(C_p):=\set{\phi_0,\phi_1,\dots,\phi_{p-1}}$, where $\phi_0 = \triv_{C_p}$. For $n=p^k$ where $k\in\N_0$, it follows from above that 
	\[ \Lin(P_{p^k}) = \set{\Chi(\psi;\phi_i) \mid \psi\in\Lin(P_{p^{k-1}}),\ i\in\set{0,1,\dotsc,p-1}}. \] 

\subsection{Plethysm coefficients}
	The irreducible characters of $S_m\wr S_n$ of the form $\Chi(\mu;\nu)$ where $\mu\in\Part(m)$ and $\nu\in\Part(n)$ are closely connected to the plethysm of Schur functions; see \cite[\textsection 7]{Stanley} or \cite{Mac95} for further detail. Let $s_\lambda$ denote the Schur function corresponding to the partition $\lambda$. The Frobenius characteristic map between class functions of finite symmetric groups and the ring of symmetric functions gives a correspondence $\chi^\lambda \longleftrightarrow s_\lambda$, and
	\[ \Chi(\mu;\nu)\ind_{S_m\wr S_n}^{S_{mn}} \longleftrightarrow s_\nu\circ s_\mu \]
	for all partitions $\lambda,\mu$ and $\nu$, where $\circ$ denotes the plethystic product. For $\lambda \in\Part(mn)$, 
	\[ \chi^\lambda\res_{S_m\wr S_n}^{S_{mn}} = \sum_{\substack{\mu\in\Part(m),\\\nu\in\Part(n)}} a^\lambda_{\nu,\mu} \Chi(\mu;\nu) \]
	where $a^\lambda_{\nu,\mu}$ is the \emph{plethysm coefficient} given by $\inprod{s_\nu\circ s_\mu, s_\lambda}$.

	\begin{lemma}[{\cite[Example 1.~(a), \textsection I.8]{Mac95}}]\label{lem: plet sgn}
		Let $\lambda$, $\mu$ and $\nu$ be partitions. Then
		\begin{align*}
			a^\lambda_{\nu,\mu} = a^{\lambda'}_{\nu^\triangle,\mu'}, \hspace{3mm}\text{ where } \nu^\triangle:=\begin{cases}\nu &\text{if } \abs{\mu} \text{ is even,} \\ \nu' &\text{if } \abs{\mu} \text{ is odd.}
		  \end{cases}
		\end{align*}
	\end{lemma}
	
	The following result is one part of \cite[Theorem 1.2]{dBPW21}.
	\begin{theorem}\label{thm:dBPW1.2}
		Let $m,n,r\in\N$ and $\mu\in\Part(m),\nu\in\Part(n),\lambda\in\Part(mn)$. Then
		\[ \inprod{s_\nu\circ s_{\mu+(1^r)}, s_{\lambda+(n^r)}} \ge \inprod{s_\nu\circ s_\mu,s_\lambda}, \]
		where partitions are added part-wise.
	\end{theorem}
	
	We also record a useful relationship between plethysm and Littlewood--Richardson coefficients: for any partition $\mu$, since $s_\mu s_\mu = s_{(2)}\circ s_\mu + s_{(1^2)}\circ s_\mu$ (see \cite{CL}, for example), then
	\begin{equation}\label{eqn:a-c}
		c^\lambda_{\mu,\mu} = a^\lambda_{(2),\mu} + a^\lambda_{(1^2),\mu}
	\end{equation}
	for all $\lambda\in\Part(2|\mu|)$.
	As will become clear in \cref{sec:ah}, we will need to compute certain values of $a^\lambda_{\nu,\mu}$ to determine the linear constituents of $\chi^\lambda\res_{P_{p^k}}$. For this, we will need the following technical lemma. 

	\begin{lemma}\label{compute a}
		Let $k\in\N_{\geq 2}$. Let $x\in\set{0,\dots,2^k-4}$, $i\in\set{0,1}$ and $\mu\in\Part(2^{k-1})$. Then
			$$ a^{\ahook{2^k}{x}}_{\hook{2}{i},\mu} = \begin{cases}
				1 &\text{if \emph{either}}\ \mu\in\set{\hook{2^{k-1}}{\frac{x}{2}}, \hook{2^{k-1}}{\frac{x+2}{2}}}\ \text{and}\ \frac{x}{2}\equiv i\ (\text{mod}\ 2)\text{, \emph{or}}\ \mu=\hook{2^{k-1}}{\frac{x+1}{2}},\\
				0 &\text{otherwise}.
			\end{cases} $$
	\end{lemma}
	\begin{proof}
		Let $\lambda=\ahook{2^k}{x}$, $\mu\in\Part(2^{k-1})$ and $W=S_{2^{k-1}}\wr S_2$. By \eqref{eqn:a-c}, we have that
      	\[ c^\lambda_{\mu,\mu} = a^\lambda_{\hook{2}{0},\mu} + a^\lambda_{\hook{2}{1},\mu} = \inprod{\lambda\res_{W}, \Chi(\mu;\hook{2}{0}) + \Chi(\mu;\hook{2}{1})}. \]
		If $\mu\not\in\set{\hook{2^{k-1}}{y}\mid y\in\set{\frac{x}{2},\frac{x+1}{2},\frac{x+2}{2}}}$, then $c^\lambda_{\mu,\mu}=0$ by the Littlewood--Richardson rule, and so $a^\lambda_{\hook{2}{0},\mu} = a^\lambda_{\hook{2}{1},\mu} = 0$. 
		
		From now on, let $\mu=\hook{2^{k-1}}{y}$ where $y\in\set{\frac{x}{2},\frac{x+1}{2},\frac{x+2}{2}}$. We first consider the case when $y=\frac{x+1}{2}$, so $x$ is odd and $c^\lambda_{\mu,\mu}=2$ by the Littlewood--Richardson rule. By repeatedly using \cref{thm:dBPW1.2} with $r=1$, we get
			$$ a^\lambda_{\hook{2}{i},\mu} \geq a^{\ahook{2^k-2}{x}}_{\hook{2}{i},\hook{2^{k-1}-1}{y}} \geq a^{\ahook{2^k-4}{x}}_{\hook{2}{i},\hook{2^{k-1}-2}{y}} \geq \dots \geq a^{\ahook{2^k-2(2^{k-1}-\frac{x+5}{2})}{x}}_{\hook{2}{i},\hook{2^{k-1}-(2^{k-1}-\frac{x+5}{2})}{y}} = a^{\ahook{x+5}{x}}_{\hook{2}{i},\hook{\frac{x+5}{2}}{y}}$$
		By \cref{lem: plet sgn} we get that $a^{\ahook{x+5}{x}}_{\hook{2}{i},\hook{\frac{x+5}{2}}{y}} = a^{\ahook{x+5}{1}}_{\hook{2}{i^\triangle},\hook{\frac{x+5}{2}}{1}}$, where $i^\triangle=i$ if $\frac{x+5}{2}$ is even, and $i^\triangle=1-i$ otherwise. Once again using \cref{thm:dBPW1.2} repeatedly with $r=1$, gives
			$$a^{\ahook{x+5}{1}}_{\hook{2}{i^*},\hook{\frac{x+5}{2}}{1}} \geq a^{\ahook{x+3}{1}}_{\hook{2}{i^*},\hook{\frac{x+3}{2}}{1}} \geq a^{\ahook{x+1}{1}}_{\hook{2}{i^*},\hook{\frac{x+1}{2}}{1}} \geq \dots \geq  a^{\ahook{x+5-(x-1)}{1}}_{\hook{2}{i^*},\hook{\frac{x+5-(x-1)}{2}}{1}} = a^{\ahook{6}{1}}_{\hook{2}{i^*},\hook{3}{1}}$$
		Since $\ahook{6}{1}$ and $\hook{3}{1}$ are self-conjugate, we have that $a^{\ahook{6}{1}}_{\hook{2}{i^*},\hook{3}{1}} = a^{\ahook{6}{1}}_{\hook{2}{i},\hook{3}{1}}$ by \cref{lem: plet sgn}. By direct computation we get $a^{\ahook{6}{1}}_{\hook{2}{i},\hook{3}{1}}=1$ for all $i\in\set{0,1}$. Hence $a^\lambda_{\hook{2}{i},\mu} \geq 1$ for all $i\in\set{0,1}$. Since $c^\lambda_{\mu,\mu}=2$, it follows that $a^\lambda_{\hook{2}{i},\mu}=1$ for all $i\in\set{0,1}$. 

    Next suppose $y=\frac{x}{2}$, so $c^\lambda_{\mu,\mu}=1$. We have that $a^\lambda_{\hook{2}{i},\mu} = a^{\ahook{2^k}{2^k-4-x}}_{\hook{2}{i},\hook{2^{k-1}}{2^{k-1}-1-y}}$ by \cref{lem: plet sgn}. Note that the first part of $\ahook{2^k}{2^k-4-x}$ is twice the size of the first part of $\hook{2^{k-1}}{2^{k-1}-1-y}$. Hence,
      $$ a^{\ahook{2^k}{2^k-4-x}}_{\hook{2}{i},\hook{2^{k-1}}{2^{k-1}-1-y}} = a^{\hook{2^k-2-x}{2^k-4-x}}_{\hook{2}{i},\hook{2^{k-1}-1-y}{2^{k-1}-2-y}}$$
    by \cite[Theorem 1.1]{dBPW21}. Now the statement follows from \cite[Theorem 3.1]{LR04}. The case when $y=\frac{x+2}{2}$ is done in a similar way. \qedhere
	\end{proof}

\section{Sylow branching coefficients for almost hooks}\label{sec:ah}
\begin{figure}[H]
	\centering
	\begin{tikzpicture}[scale=0.8, every node/.style={scale=0.75}]
		\draw (0,0) node(S) {$S_{2^k}$};
		\draw (0,-1) node(W) {$W=S_{2^{k-1}}\wr S_2$};
		\draw (-3,-2) node(Y) {$Y=S_{2^{k-1}}\times S_{2^{k-1}}$};
		\draw (3,-2) node(P) {$P=P_{2^k}=P_{2^{k-1}}\wr P_2 = Q \wr P_2$};
		\draw (0,-3) node(B) {$B=P_{2^{k-1}}\times P_{2^{k-1}} = Q \times Q$};
		\draw (S) -- (W) -- (Y) -- (B) -- (P) -- (W);
	\end{tikzpicture}
\caption{Some subgroups of $S_{2^k}$.}\label{fig: diamond}
\end{figure}
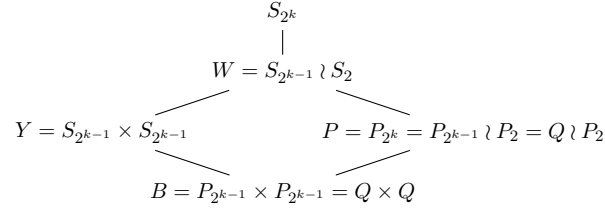

	The main aim of this section is to prove \cref{thm:almost-hook-sbcs}. We fix $p=2$ in this section.
	We begin with a simple but important observation. Let $k\in\N$, $K\leqslant H \leqslant S_{2^k}$ and $\chi\in\Irr(S_{2^k})$. Firstly, if $\psi\in\Irr(K)$ is a constituent of $\chi\res_K$, then there must exist some $\theta\in\Irr(H)$ such that $\theta$ is a constituent of $\chi\res_H$ and $\psi$ is a constituent of $\theta\res_K$. 
	The subgroups of $S_{2^k}$ that will be of interest here are those in \cref{fig: diamond}; we will use this notation throughout the rest of this paper. 

	Let $k\geq 2$ and $\lambda\in\Part(2^k)$. Suppose that there exists $\mu,\nu\in\Part(2^{k-1})$ such that $c^\lambda_{\mu,\nu}>0$ and $\phi\in\Lin(\chi^\mu\res_Q)\cap\Lin(\chi^\nu\res_Q)$ for some $\phi\in\Lin(Q)$, then $\inprod{\chi^\lambda\res_{B},\phi\times\phi}>0$. Consequently, there must exist $\theta\in\Irr(P)$ such that $\inprod{\theta\res_B,\phi\times\phi}>0$ and $\inprod{\chi^\lambda\res_P,\theta}>0$. Thus $\theta\in \Irr(P\mid\phi\times\phi) = \set{\Chi(\phi;\phi_0),\Chi(\phi;\phi_1)}$.
	
	Let $Y(\lambda)$ denote the set of $\phi\in\Lin(P_{2^{k-1}})$ such that $c^\lambda_{\mu,\nu}>0$ and $\phi\in\Lin(\chi^\mu\res_Q)\cap\Lin(\chi^\nu\res_Q)$ for some $\mu,\nu\in\Part(2^{k-1})$. Then
		\begin{align}\label{eq:find lins}
			\abs{\Lin(\chi^\lambda\res_P)} \geq \abs{Y(\lambda)}
		\end{align}
	since $\Chi(\phi;\phi_0)\res_B = \Chi(\phi;\phi_1)\res_B = \phi\times\phi$. In order to determine for which $i\in\set{0,1}$ we have that $\Chi(\phi;\phi_i)\in\Lin(\chi^\lambda\res_P)$ we will look at constituents of $\chi^\lambda\res_W$, which is done by considering relevant plethysm coefficients as introduced in \cref{sec:prelims}. 
	
	Next, we record a useful bijection from \cite[Theorem 1.1]{G} and \cite[Theorem 3.1]{GV-vietnam}, stated in present notation. 

	\begin{theorem}\label{thm: hook lins}
		Let $k\in \N$ and $\lambda\in\Hook(2^k)$. Then $\chi^\lambda\res_{P_{2^k}}$ has a unique linear constituent $\lin_{k}(\lambda)$. Moreover, the map
			$$ \lin_k: \Hook(2^k) \rightarrow \Lin(P_{2^k})$$
		is a bijection, and $\inprod{\chi^\lambda\res_{P_{2^k}}, \lin_k(\lambda)}=1$ for all $\lambda\in\Hook(2^k)$.
	\end{theorem}
	When clear from context, we will abbreviate $\lin_k$ to $\lin$.\footnote{\cref{thm: hook lins} in fact describes an instance of a McKay bijection $\Irr_{p'}(G)\to\Irr_{p'}(N_G(P))$, where $G=S_{2^k}$ and $p=2$. This is because Sylow 2-subgroups of symmetric groups are self-normalising, $\Lin(P_{2^k})=\Irr_{2'}(P_{2^k})$, and $\Irr_{2'}(S_{2^k}) = \{\chi^\lambda\mid\lambda\in\Hook(2^k)\}$; further details appear in \cite{OlssonBook}.} The following is an immediate consequence of \cite[Theorem 3.2]{GV-vietnam}\label{thm: giavol}.

	\begin{proposition}\label{rewrite lin}
		Let $k\in\N$, $i\in\set{0,\dots,2^{k}-1}$ and $j\in\{0,1\}$. Then
		\[ \Chi( \lin(\hook{2^k}{i}); \phi_j ) = \begin{cases}
			\lin(\hook{2^{k+1}}{2i+j}) & \text{if $i$ is even},\\
			\lin(\hook{2^{k+1}}{2i+1-j}) & \text{if $i$ is odd}.
		\end{cases} \]
	\end{proposition}

	We are now ready to prove \cref{thm:almost-hook-sbcs}. For $a\in\N_0$ and $b\in\Z$, recall ${a\choose b}=0$ if $b\not\in\set{0,\dots,a}$.

	\begin{proof}[Proof of \cref{thm:almost-hook-sbcs}]
		We proceed by induction on $k$. For $k\le 3$, the statement follows from direct computation. Now assume $k\ge 3$ and that the statement holds true for $k-1$. Let $\lambda=\ahook{2^k}{x}$ and $\gamma = \hook{2^k}{y}$. We need to compute $\sbc^\lambda_{\lin(\gamma)}$.
		Following the notation of \cref{fig: diamond}, recall the description of $\Irr(W)$ in \cref{sec:prelims}. Restriction of $\chi^\lambda$ to $P$ through $W$ gives us
		\begin{align}
			\sbc^{\lambda}_{\lin(\gamma)} &= \sum_{\substack{\mu\in\Part(2^{k-1}),\\ i\in\set{0,1}}} \inprod{\lambda\res_W, \Chi(\mu; \phi_i)}\cdot\inprod{\Chi(\mu;\phi_i)\res_P, \lin(\gamma)} \label{first inprod} \\
			&\phantom{=0} + \frac{1}{2}\sum_{\substack{\eta,\nu\in\Part(2^{k-1}),\\ \eta\neq \nu}} \inprod{\lambda\res_W, (\eta\times\nu)\ind^W}\cdot\inprod{(\eta\times\nu)\ind^W\res_P, \lin(\gamma)} \label{second inprod}
		\end{align}
		(Note the factor of $\tfrac{1}{2}$ in \eqref{second inprod} arises since $(\eta\times\nu)\ind_Y^W = (\nu\times\eta)\ind_Y^W$.) 
		We proceed to consider the terms in \eqref{first inprod} and \eqref{second inprod} separately. 
		Using \cref{rewrite lin}, let $j\in\set{0,1}$ be such that $\lin(\gamma)=\Chi(\lin(\hook{2^{k-1}}{z});\phi_j)$, where $z=\lfloor\frac{y}{2}\rfloor$. 
		\begin{itemize}
			\item For \eqref{first inprod} we note that $\inprod{\lambda\res_W, \Chi(\mu;\phi_i)} = a^\lambda_{\hook{2}{i},\mu}$ and $a^\lambda_{\hook{2}{i},\mu}=0$ for all $\mu\in\Part(2^{k-1})\setminus\set{\hook{2^{k-1}}{\frac{x+t}{2}}\mid t\in\set{0,1,2}}$ by \cref{compute a}.
			Moreover, for $\mu\in\Hook(2^{k-1})$, 
			\[ \inprod{\Chi(\mu;\phi_i)\res_P, \lin(\gamma)} =\inprod{\Chi(\mu\res_Q;\phi_i), \Chi(\lin(\hook{2^{k-1}}{z});\phi_j)} = \delta_{\mu,\hook{2^{k-1}}{z}} \cdot \delta_{i,j},\]
			where the final equality follows from \cref{thm: hook lins} and \cite[Lemma 2.19]{SLThesis}.
			
			\item Next consider \eqref{second inprod}. Firstly, note that $\inprod{\lambda\res_W,(\eta\times\nu)\ind^W} = c^\lambda_{\eta,\nu}$. Since $\lambda$ is an almost hook, if $c^\lambda_{\eta,\nu} > 0$ then we may without loss of generality assume that $\nu\in\Hook(2^{k-1})$ and $\eta\in\Hook(2^{k-1})\cup \Ahook(2^{k-1})$. By Mackey's Theorem (see \cite[\textsection 5]{Isaacs}), we have that $(\eta\times\nu)\ind^W\res_P = (\eta\times\nu)\res_B\ind^P$, and so
			\begin{align*}
				\inprod{(\eta\times\nu)\ind^W\res_P, \lin(\gamma)} &= \inprod{(\eta\times\nu)\res_B\ind^P,\lin(\gamma)} = \inprod{(\eta\times\nu)\res_B,\lin(\gamma)\res_B} \\
				&= \inprod{\eta\res_Q \times \nu\res_Q, \lin(\hook{2^{k-1}}{z})\times \lin(\hook{2^{k-1}}{z})} \\
				&= \sbc^\eta_{\lin(\hook{2^{k-1}}{z})} \cdot \delta_{\nu,\hook{2^{k-1}}{z}},
			\end{align*}
			where the final equality follows from \cref{thm: hook lins} since $\nu\in\Hook(2^{k-1})$. Moreover, $\eta\ne\nu$, so if $\eta\in\Hook(2^{k-1})$ then $\sbc^\eta_{\lin(\hook{2^{k-1}}{z})}=0$ by \cref{thm: hook lins}. Thus we have non-zero contributions to \eqref{second inprod} only if, without loss of generality, $\nu\in\Hook(2^{k-1})$ and $\eta\in\Ahook(2^{k-1})$.
		\end{itemize}
		Simplifying the expression for $\sbc^\lambda_{\lin(\gamma)}$ in terms of \eqref{first inprod} and \eqref{second inprod}, we obtain
		\[ \sbc^\lambda_{\lin(\gamma)} = a^\lambda_{\hook{2}{j},\hook{2^{k-1}}{z}} + \sum_{\eta\in\Ahook(2^{k-1})} c^\lambda_{\eta,\hook{2^{k-1}}{z}}\cdot \sbc^\eta_{\lin(\hook{2^{k-1}}{z})}. \]
		Therefore, by the Littlewood--Richardson rule,
		\begin{equation}\label{eqn:Z}
			\sbc^\lambda_{\lin(\gamma)} = a^\lambda_{\hook{2}{j},\hook{2^{k-1}}{z}} + \delta_0\cdot \sbc^{\ahook{2^{k-1}}{x-z}}_{\lin(\hook{2^{k-1}}{z})} + \delta_1\cdot \sbc^{\ahook{2^{k-1}}{x-z-1}}_{\lin(\hook{2^{k-1}}{z})}
		\end{equation}
		where for each $l\in\set{0,1}$, $\delta_l = 1$ if $x-z-l\in\set{0,1,\dotsc,2^{k-1}-4}$ and $\delta_l=0$ otherwise.
		Note that since $y\in\set{2z,2z+1}$, we have $b(z)-2z=b(y)-y$ and so $B(z)+z = B(y)$. Substituting \cref{compute a} and the induction hypothesis into \eqref{eqn:Z}, we find that  
		\begin{align*}
			\sbc^\lambda_{\lin(\gamma)} &= \begin{cases}
				1 &\text{if \emph{either}}\ x\in\set{2z,2z-2} \text{ and }\frac{x}{2}\equiv j\ (\text{mod}\ 2) \text{, \emph{or} } x = 2z-1,\\
				0 &\text{otherwise};
			\end{cases}\\
			&+ \begin{cases}
				{k-2 \choose x-B(y)}-1 &\text{if } x\in\set{2z-1, 2z-2},\\
				{k-2 \choose x-B(y)} &\text{if } x\in\set{B(y),\dots,B(y)+k-2}\setminus\set{2z-1,2z-2},\\
				0 &\text{otherwise};
			\end{cases}\\
			&+ \begin{cases}
				{k-2 \choose x-B(y)-1}-1 &\text{if } x\in\set{2z, 2z-1},\\
				{k-2 \choose x-B(y)-1} &\text{if } x\in\set{B(y)+1,\dots,B(y)+k-1}\setminus\set{2z,2z-1},\\
				0 &\text{otherwise.}
			\end{cases}
		\end{align*}
		Finally, we analyse the expression above by considering the cases $x=2z-1$, $x\in\set{2z,2z-2}$ and $x\in\set{B(y),\dots,B(y)+k-1}\setminus\set{2z,2z-1,2z-2}$ separately. 
		Notice since $y\in\set{0,1,\dotsc,2^k-1}$ then $b(y)\le k$ and so $B(y)\in\set{y-k-1,\dots,y-1}$.
		If $x\in\set{B(y),\dots,B(y)+k-1}\setminus\set{2z,2z-2}$, then it is straightforward to verify that the statement is true, since $B(y)\le 2z-1\le B(y)+k-1$ and $\binom{N}{t}=\binom{N-1}{t}+\binom{N-1}{t-1}$. 
		If $x=y=2z$, then $\frac{x}{2} = z$ and so $\frac{x}{2}\equiv j$ (mod $2$) by \cref{rewrite lin}. It follows that
		\begin{align*}
			\sbc^\lambda_{\lin(\gamma)}= 1 + {k-2 \choose x-B(y)} + {k-2 \choose x-B(y)-1} -1 = {k-1 \choose x-B(y)}
		\end{align*}
		as desired. If $x=2z$ and $y=2z+1$ then $\frac{x}{2}=z\not\equiv j$ (mod $2$) and $\sbc^\lambda_{\lin(\gamma)}=\binom{k-1}{x-B(y)}-1$ as desired. The case $x=2z-2$ is similar.
	\end{proof}

We can now use \cref{thm:almost-hook-sbcs} to compute $\abs{\Lin(\chi^\lambda\res_P)}$ for $\lambda\in\Ahook(2^k)$, where $P\in\Syl_2(S_{2^k})$.

\begin{corollary}\label{cor: ah lins}
	Let $k\in\N_{\geq 2}$, $P\in\Syl_2(S_{2^k})$ and $x\in\set{0,\dots,2^k-4}$.
	\begin{enumerate}[label=\roman*)]
		\item If $x = 2^{k-1}-2$, then $\Lin(\ahook{2^k}{x}\res_P) = \set{ \lin(\hook{2^k}{2^{k-1}-2}), \lin(\hook{2^k}{2^{k-1}+1}) }$.
		Moreover, $\inprod{\ahook{2^k}{x},\psi}=1$ for all $\psi\in\Lin(\ahook{2^k}{x}\res_P)$.
		\item If $x\ne 2^{k-1}-2$, then $\abs{\Lin(\ahook{2^k}{x}\res_P)} > 2$.
	\end{enumerate}
\end{corollary}

\begin{proof}
	We abbreviate $\sbc^{\ahook{2^k}{x}}_{\lin(\hook{2^k}{y})}=:\sbc^x_y$. Note \eqref{eqn:sgn} gives $\lin(\hook{2^k}{2^k-1-y})=\sgn_P\cdot\lin(\hook{2^k}{y})$, and so $\sbc^x_y = \sbc^{2^k-4-x}_{2^k-1-y}$. Moreover, $z\ge b(z)$ for all $z\in\N_0$, with equality if and only if $z\in\{0,1\}$.
	
	\noindent \textbf{i)} First suppose $x=2^{k-1}-2$. It is easy to check that $\sbc^x_y>0$ when $y=x$, thus $\sbc^x_y>0$ also for $y=2^{k-1}+1$. Hence $|\Lin(\ahook{2^k}{x}\res_P)|\ge 2$ and moreover $\sbc^x_y=1$ for $y\in\{2^{k-1}-2,2^{k-1}+1\}$. 
	
	Conversely, suppose $\sbc^x_y>0$. We may assume without loss of generality that $y\le 2^{k-1}-1$. By \cref{thm:almost-hook-sbcs}, $x\in\{B(y),\dotsc,B(y)+k-1\}$, or equivalently $2^{k-1}-k\le y-b(y)\le 2^{k-1}-1$. In particular, $y\ge 2^{k-1}-k$.
	If $y=2^{k-1}-1=x+1$ then $\sbc^x_y=0$, and we have already considered $y=2^{k-1}-2$ above, so we can now suppose $y=2^{k-1}-1-z$ where $2\le z\le k$. Notice $b(y)=k-1-b(z)$, so
	\[ y-b(y)=2^{k-1}-1-z-(k-1)+b(z)=2^{k-1}-k - (z-b(z))<2^{k-1}-k, \]  a contradiction. Hence $\Lin(\ahook{2^k}{x}\res_P) = \set{ \lin(\hook{2^k}{2^{k-1}-2}), \lin(\hook{2^k}{2^{k-1}+1}) }$.
	
	\noindent \textbf{ii)} Now suppose $x\ne2^{k-1}-2$. It suffices to show for each $0\le x<2^{k-1}-2$ that $\sbc^x_y>0$ for three different values of $y\in\{0,\dotsc,2^k-1\}$. If $x=2^{l-1}-2$ for some $2\le l<k$, it is straightforward to check that $\sbc^x_y>0$ for $y\in\{x,x+1,x+3\}$. Similarly, if $x\ne 2^{l-1}-2$ for any $2\le l<k$ then $\sbc^x_y>0$ for $y\in\{x,x+1,x+2\}$.
\end{proof}

\begin{remark}\label{rem:identify-two-linears}
	Let $L:=\Lin(\chi^\lambda\res_P)$. As previously noted, when $p=2$ it was determined in \cite[Theorem 1.2]{G} when $\abs{L}=1$: see \cref{thm: hook lins} above. To identify the elements of $L$ when $\abs{L}=2$: this is given by \cref{cor: ah lins} in case (1) of \cref{thm:main-p=2}. In case (2), $\Lin(\hook{2^k+1}{t}\down_{P_{2^k+1}}) = \set{ \lin(\hook{2^k}{t})\times\triv_{P_1}, \lin(\hook{2^k}{t-1})\times\triv_{P_1} }$, from \cite[Remark 3.13]{GN}. In case (3), we find by direct computation
	\begin{small}
		\[ \begin{array}{cccccl}
			\Lin((5,3)\down_{P_8}) &=& \{ \lin(\hook{8}{1}), \lin(\hook{8}{2}) \}, & \Lin((3^2,2)\down_{P_8}) &=& \{ \lin(\hook{8}{2}), \lin(\hook{8}{5}) \}, \\
			\Lin((2^3,1^2)\down_{P_8}) &=& \{ \lin(\hook{8}{5}), \lin(\hook{8}{6}) \}, & \Lin((3^3)\down_{P_9}) &=& \{ \lin(\hook{8}{2})\times \triv_{P_1}, \lin(\hook{8}{5})\times \triv_{P_1} \}.
		\end{array} \]
	\end{small}
	
\end{remark}

\section{Proof of \cref{thm:main-p=2}}\label{sec:p=2}
In this section, fix $p=2$. We determine for which $\chi\in\Irr(S_n)$ we have that $\abs{\Lin(\chi\res_{P_n})}=2$. We begin with a technical lemma that will be useful for finding a lower bound on $\abs{\Lin(\chi\res_{P_n})}$. Given two partitions $\nu=(\nu_1,\nu_2,\dotsc)$ and $\lambda=(\lambda_1,\lambda_2,\dotsc)$, we say $\nu\subseteq\lambda$ if $\nu_i\le\lambda_i$ for all $i$.

\begin{lemma}\label{lem: lins lower bound}
	Let $n\in\N$, $p$ be a prime and let $n=a_0p^0 + \dots + a_kp^k$ be the $p$-adic expansion of $n$. Let $m\in\N$ be such that $m$ has $p$-adic expansion $m=b_0p^0+\cdots+b_lp^l$ with $l\le k$ and $b_i\le a_i$ for all $i$. If $\nu\in\Part(m)$ and $\lambda\in\Part(n)$ satisfy $\nu\subseteq \lambda$, then $\abs{\Lin(\chi^\lambda\res_{P_n})}\geq \abs{\Lin(\chi^\nu\res_{P_m})}$.
\end{lemma}

\begin{proof}
	Since $\nu\subseteq \lambda$, there exists some $\gamma\in\Part(n-m)$ such that $c^\lambda_{\nu,\gamma}>0$. Hence we have that $\Lin((\nu\times \gamma)\res_{P_n})\subseteq \Lin(\lambda\res_{P_n})$, as $P_n\cong P_m\times P_{n-m}\le S_m\times S_{n-m}$ by \eqref{eqn:Pn}. The claim follows since $\abs{\Lin(\gamma\res_{P_{n-m}})}\ge1$ by \cref{prop: at least one lin} and $\abs{\Lin((\nu\times\gamma)\res_{P_n})}=\abs{\Lin(\nu\res_{P_m})}\cdot\abs{\Lin( \gamma\res_{P_{n-m}})}$.
\end{proof}

Let $m\in\N$, $\mu\in\Part(m)$ and $\lambda\in\Part(n)$. We denote by $\mu\sqcup\lambda$ the partition of size $m+n$ whose multiset of parts is the union of those of $\mu$ and $\lambda$. Furthermore, if $\gamma$ is a composition of $m$ (i.e.~a sequence of non-negative integers whose sum is $m$), then let $\gamma^*$ denote the partition obtained by rearranging the parts of $\gamma$, omitting any zeros if necessary. For example, if $\lambda=(5,3,1)\in\Part(9)$ and $\gamma=(0,3,2,7)$ then $\gamma^*=(7,3,2)\in\Part(12)$ and $\lambda\sqcup\gamma^* = (7,5,3,3,2,1)\in\Part(21)$. 

\begin{definition}\label{delta part}
	Let $n\in \N$ and let $\lambda=(\lambda_1,\dots,\lambda_{\ell(\lambda)})\in\Part(2n)$.
		\begin{itemize}
			\item Let $\alpha=(\alpha_1,\dots,\alpha_{\ell(\alpha)})$ be the partition consisting of the parts of $\lambda$ which are even.
			\item Let $\beta=(\beta_1,\dots,\beta_{\ell(\beta)})$ be the partition consisting of the parts of $\lambda$ which are odd. Notice that $\ell(\beta)$ must be even since $\lambda\in\Part(2n)$; say $\ell(\beta)=2t$ for some $t\in\N_0$.
		\end{itemize}
	We then define the partition $\Delta(\lambda)\in\Part(n)$ by
	\begin{align*} \Delta(\lambda) = \frac{\alpha}{2} \sqcup \left(\frac{\beta_1+1}{2}, \frac{\beta_2-1}{2}, \frac{\beta_3+1}{2}, \frac{\beta_4-1}{2}, \cdots, \frac{\beta_{2t-1}+1}{2}, \frac{\beta_{2t}-1}{2}\right)^*. \end{align*}
	We remark that this coincides with \cite[Definition 4.1]{G}.
\end{definition}

\begin{proposition}\cite[Proposition 4.3]{G}\label{split part}
	Let $n\in\N$ and $\lambda\in\Part(2n)$, then $c^\lambda_{\Delta(\lambda),\Delta(\lambda)}>0$.
\end{proposition}

Together with \eqref{eq:find lins}, \cref{split part} gives a lower bound 
\[ \abs{\Lin(\chi^\lambda\res_{P_{2^k}})} \geq \abs{\Lin(\chi^{\Delta(\lambda)}\res_{P_{2^{k-1}}})}.\]
If $\alpha\in\Hook(2^{k-1})\cup\set{\ahook{2^{k-1}}{2^{k-2}-2}}$ then $\abs{\Lin(\chi^\alpha\res_{P_{2^{k-1}}})}\le 2$ by \cref{thm: hook lins} and \cref{cor: ah lins}. 
To prove \cref{thm:main-p=2} when $n=2^k$ by induction on $k$, the main idea is that for all other $\alpha\in\Part(2^{k-1})$, with a small number of exceptions, we will have $\abs{\Lin(\chi^\alpha\res_{P_{2^{k-1}}})}>3$ by the inductive hypothesis. Moreover, recall that $\abs{\Lin(\chi^\lambda\res_{P_n})} = \abs{\Lin(\chi^{\lambda'}\res_{P_n})}$. Hence, we want to determine those $\lambda\in\Part(2^k)$ for which $\set{\Delta(\lambda),\Delta(\lambda')}\subseteq\Hook(2^{k-1})\cup\set{\ahook{2^{k-1}}{2^{k-2}-2}}$. 

\begin{definition}\label{def:Dk}
	Let $k\in\N_{\ge 3}$. We define $\cD_k$ to be the subset of $\Part(2^k)$ consisting of those partitions of the form $\mu\sqcup\nu\sqcup(1^{2^k-\vert\mu\vert-\vert\nu\vert})$, where $\mu_1\ge\nu_1$ and $\mu$ and $\nu$ satisfy one of the following:
	\begin{enumerate}[label=\alph*)]
		\item $\mu=(2^{k-1})$ and $\nu\in\set{(3), (4),(4,2)}$,
		\item $\mu=(2^{k-1}-1)$ and $\nu\in\set{(4),(5),(5,2),(2^2),(2^3)}$,
		\item $\mu=(2^{k-1}-2)$ and $\nu\in\set{(2^3),(2^4)}$,
		\item $\mu=(2^{k-1}-3)$ and $\nu\in\set{(3,2^2),(3,2^3)}$,
		\item $\mu=(2r+3)$ for some $r\in\N_0$ and $\nu\in\set{(3),(3,2)}$,
		\item $\mu=(2r)$ for some $r\in\N$ and $\nu=(2^2)$.
	\end{enumerate} 
\end{definition}

\begin{lemma}\label{annoying partitions}
	Let $k\in\N_{\ge 3}$. Then $\lambda\in\cD_k\cup\Hook(2^k)\cup\Ahook(2^k)$ if and only if $\set{\Delta(\lambda),\Delta(\lambda')}\subseteq\Hook(2^{k-1})\cup\set{\ahook{2^{k-1}}{2^{k-2}-2}}$. 
\end{lemma}

\begin{proof}
	Suppose that $\lambda\in\cD_k\cup\Hook(2^k)\cup\Ahook(2^k)$. It is straightforward to verify that $\Delta(\lambda),\Delta(\lambda')\in\Hook(2^{k-1})\cup\set{\ahook{2^{k-1}}{2^{k-2}-2}}$. 
	For the converse, we first suppose $\lambda\in\Part(2^k)$ is such that $\Delta(\lambda)=\ahook{2^{k-1}}{2^{k-2}-2}$. 
	Then $\Delta(\lambda)_1=2^{k-2}$, which implies that $\lambda_1\in\set{2^{k-1},2^{k-1}-1}$. Furthermore, $\Delta(\lambda)_2 = 2$ implies that $\lambda_2=4$, \textit{or} $\lambda_2=3$ and $\lambda_1$ is even, \textit{or} $\lambda_2=5$ and $\lambda_1$ is odd. Lastly, we have that $\Delta(\lambda)_i=1$ for $i\in\set{3,\dots,\ell(\Delta(\lambda))}$. It follows that $\lambda_i\in\set{1,2,3}$ for $i\in\set{3,\dots,\ell(\lambda)}$, with $\lambda_i=3$ only if $\lambda_1$ and $\lambda_2$ have different parities. Hence $\lambda\in\cA$ where $\cA\subseteq\cP(2^k)$ comprises those partitions of the following forms:
	\begin{multicols}{2}
		\begin{enumerate}[label=(\alph*)]
			\item $(2^{k-1},4,2^a,1^{2^{k-1}-4-2a})$,
			\item $(2^{k-1},3^2,2^a,1^{2^{k-1}-6-2a})$,
			\item $(2^{k-1},3,2^a,1^{2^{k-1}-3-2a})$,
			\item $(2^{k-1}-1,4,2^a,1^{2^{k-1}-3-2a})$,
			\item $(2^{k-1}-1,4,3,2^a,1^{2^{k-1}-6-2a})$,
			\item $(2^{k-1}-1,5,2^a,1^{2^{k-1}-4-2a})$,
		\end{enumerate}
	\end{multicols}
	\noindent where $a\in\N_0$ (and $k\ge 4$ in forms (d), (e) and (f)). 
	Similarly, we find that if $\lambda\in\Part(2^k)$ satisfies $\Delta(\lambda)\in\Hook(2^{k-1})$ then $\lambda\in\cB\cup\Hook(2^k)$, where 
	\begin{align*}
		\cB &:= \set{(r,2^a,1^{2^k-r-2a})\in\Part(2^k)\mid r\in\N_{\ge 2}, a\in\N_0} \\
		& \qquad \cup \set{(r,3,2^a,1^{2^k-3-r-2a})\in \Part(2^k)\mid r\in\N_{\geq 3}\text{ and $r$ is odd}, a\in\N_0}.
	\end{align*}
	It is then straightforward to verify that $\set{\lambda,\lambda'}\subseteq\cA\cup\cB\cup\Hook(2^k)$ implies $\lambda\in\cD_k\cup\Hook(2^k)\cup\Ahook(2^k)$.
\end{proof}

\begin{lemma}\label{lem:fixing spc cases}
	Let $k\in\N_{\geq 4}$. If $\lambda\in\cD_k{\blue\setminus\{(2^k-3,3),(2^3,1^{2^k-6})\}}$ then there exist $\alpha,\beta\in\Hook(2^{k-1})\cup\Ahook(2^{k-1})$ such that $c^\lambda_{\alpha,\alpha},c^\lambda_{\beta,\beta}>0$ and $\abs{\Lin(\alpha\res_{P_{2^{k-1}}}) \cup \Lin(\beta\res_{P_{2^{k-1}}})} > 2$.
\end{lemma}

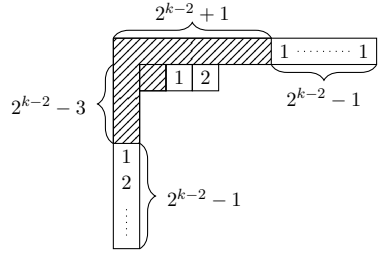
\begin{figure}[h!]
	\centering
	\begin{tikzpicture}[scale=0.35,every node/.style={scale=.7}]
		\draw (0,0) -- (10,0) -- (10,-1) -- (2,-1) -- (2,-2) -- (1,-2) -- (1,-8) -- (0,-8) -- (0,0);
		\draw [pattern=north east lines] (0,0) -- (6,0) -- (6,-1) -- (1,-1) -- (1,-2) -- (1,-4) --(0,-4) -- (0,0);
		\draw [pattern=north east lines] (1,-1) -- (1,-2) -- (2,-2) -- (2,-1) -- (1,-1);
		\draw (2,-1) -- (2,-2) -- (3,-2) -- (3,-1) -- (2,-1);
		\draw (3,-1) -- (3,-2) -- (4,-2) -- (4,-1) -- (3,-1);
		\draw [decoration={brace,amplitude=6pt},decorate] (0,0) -- (6,0) node[midway,yshift=14]{$2^{k-2}+1$};
		\draw [decoration={brace,amplitude=6pt},decorate] (0,-4) -- (0,-1) node[midway, xshift=-35,yshift=-1]{$2^{k-2}-3$};
		\draw [decoration={brace,amplitude=6pt},decorate] (10,-1) -- (6,-1) node[midway,yshift=-18]{$2^{k-2}-1$};
		\draw [decoration={brace,amplitude=6pt},decorate] (1,-4) -- (1,-8) node[midway,xshift=35,yshift=-1]{$2^{k-2}-1$};
		\draw (6.5,-0.5) node[] {$1$};
		\draw [dotted] (7,-.5) -- (9,-.5);
		\draw (9.5,-.5) node[] {$1$};
		\draw (2.5,-1.5) node[] {$1$};
		\draw (3.5,-1.5) node[] {$2$};
		\draw (.5,-4.5) node[] {$1$};
		\draw (.5,-5.5) node[] {$2$};
		\draw [dotted] (.5,-6.5) -- (.5,-7.5);
	\end{tikzpicture}
	\caption{Using the Littlewood--Richardson rule to verify $c^\lambda_{\alpha,\alpha}>0$ where $\lambda=\mu\sqcup\nu\sqcup(1^{2^k-|\mu|-|\nu|})$, $\mu=(2^{k-1})$, $\nu=(4)$ and $\alpha=\ahook{2^{k-1}}{2^{k-2}-3}$: the shape $\lambda$ is split into shaded and unshaded parts. The shaded part illustrates $\alpha$, and the numbers filled into the unshaded part have reverse row word of content $\alpha$.} \label{fig:LR-example}
\end{figure}

\begin{proof}
 	Let $Q=P_{2^{k-1}}$ and write $\lambda=\mu\sqcup\nu\sqcup(1^{2^k-|\mu|-|\nu|})$ as in \cref{def:Dk}. 
 	\begin{itemize}
 		\item For $\lambda$ in \cref{tab:alpha}, we may take $\alpha=\beta$ as listed, noting that $\abs{\Lin(\alpha\res_Q)} \ge 3$ by \cref{cor: ah lins}. We note that $c^\lambda_{\alpha,\alpha}>0$ by the Littlewood--Richardson rule: see \cref{fig:LR-example} for an example illustration in the case where $\mu=(2^{k-1})$, $\nu=(4)$ and $\alpha=\ahook{2^{k-1}}{2^{k-2}-3}$. The other cases can be verified similarly.
 		
 		\item For $\lambda$ in \cref{tab:alphabeta}, we may take $\alpha$ as listed and $\beta=\ahook{2^{k-1}}{2^{k-2}-2}$, noting that $\abs{\Lin(\alpha\res_Q)}=1$, $\abs{\Lin(\beta\res_Q)}=2$ and $\Lin(\alpha\res_Q)\cap\Lin(\beta\res_Q)=\varnothing$ by \cref{thm: hook lins,thm:almost-hook-sbcs}. Again, $c^\lambda_{\alpha,\alpha}>0$ and $c^\lambda_{\beta,\beta}>0$ follow from the Littlewood--Richardson rule.
 	\end{itemize}
\end{proof}

	{\renewcommand{\arraystretch}{1.2}
	\begin{table}[h!]
		\parbox{.58\linewidth}{
			\centering
			\begin{footnotesize}
				\begin{tabular}{ll|l}
					\hline
					$\mu$ & $\nu$ & $\alpha$ \\
					\hline
					\hline
					$(2^{k-1})$ & $(4)$, $(4,2)$ & $\ahook{2^{k-1}}{2^{k-2}-3}$ \\
					$(2^{k-1}-1)$ & $(5)$, $(5,2)$ & \\
					\hline
					$(2^{k-1}-2)$ & $(2^3)$, $(2^4)$ & $\ahook{2^{k-1}}{2^{k-2}-1}$ \\
					$(2^{k-1}-3)$ & $(3,2^2)$, $(3,2^3)$ & \\
					\hline
					$(2r+3)$, {\tiny $r\ne 2^{k-2}-2$} & $(3)$, $(3,2)$ & $\ahook{2^{k-1}}{2^{k-1}-r-4}$ \\
					\hline
					$(2r)$, {\tiny $r\ne 2^{k-2}$} & $(2^2)$ & $\ahook{2^{k-1}}{2^{k-1}-r-2}$ \\
					\hline
				\end{tabular}
				\caption{\footnotesize{$\alpha=\beta$ for some $\lambda$ in \cref{lem:fixing spc cases}. {\blue For $r=2^{k-2}-2$ in the second-last row (resp.~$r=2^{k-2}$ in the last row), the corresponding $\lambda$ is covered in the last (resp.~first) row of \cref{tab:alphabeta}.}}}
				\label{tab:alpha}
			\end{footnotesize}
		}
		\hfill
		\parbox{.4\linewidth}{
			\centering
			\begin{footnotesize}
				\begin{tabular}{ll|l}
					\hline
					$\mu$ & $\nu$ & $\alpha$ \\
					\hline
					\hline
					$(2^{k-1})$ & $(3)$, $(2^2)$ & $\hook{2^{k-1}}{2^{k-2}-1}$ \\
					$(2^{k-1}-1)$ & $(4)$ & \\
					\hline
					$(2^{k-1}-1)$ & $(2^2)$, $(2^3)$ & $\hook{2^{k-1}}{2^{k-2}}$ \\
					& $(3)$, $(3,2)$ & \\
					\hline
				\end{tabular}
				\caption{\footnotesize{$\alpha$ and $\beta=\ahook{2^{k-1}}{2^{k-2}-2}$ for some $\lambda$ in \cref{lem:fixing spc cases}}}
				\label{tab:alphabeta}
			\end{footnotesize}
		}
	\end{table}
	}

{\blue
\begin{lemma}\label{lem:2k-3}
	Let $k\in\N_{\ge 3}$ and $x\in\{0,1,\dotsc,2^k-1\}$. Then
	\[ \sbc^{(2^k-3,3)}_{\cL(\hook{2^k}{x})} = \begin{cases}
		\tfrac{k(k-3)}{2}+x & \text{if }x\in\{0,1\},\\
		k-x & \text{if }x\in\{2,3\},\\
		0 & \text{otherwise}.
	\end{cases} \]
\end{lemma}
\begin{proof}
	We proceed by induction on $k$, noting that the case $k=3$ holds by direct computation. Recall the subgroups of $S_{2^k}$ from \cref{fig: diamond}, and note $\phi_r=\chi^{\hook{2}{r}}\in\Irr(C_2)$ for $r\in\{0,1\}$. Let $\lambda=(2^k-3,3)$ and $\sL=\cL(\hook{2^k}{x})$. By \cref{rewrite lin}, $\sL=\Chi\big(\cL(\hook{2^{k-1}}{z});\phi_i\big)$ for some $i\in\{0,1\}$ where $z:=\lfloor\tfrac{x}{2}\rfloor$. Let $>$ denote the lexicographic order on $\cP(2^{k-1})$. Then, similarly to the proof of \cref{thm:almost-hook-sbcs}, we have that
	\begin{footnotesize}
	\begin{align*}
		\sbc^{\lambda}_{\sL} &= \sum_{\psi\in\Irr(W)} \inprod{\lambda\res_W, \psi} \cdot \inprod{\psi\res_P, \sL}\\
		&= \sum_{\substack{\mu\in\cP(2^{k-1}),\\ r\in\{0,1\}}} \inprod{\lambda\res_W, \Chi(\mu;\phi_r)} \cdot \inprod{\Chi(\mu;\phi_r)\res_P, \sL} + \sum_{\substack{\eta,\nu\in\cP(2^{k-1}),\\ \eta>\nu}} \inprod{\lambda\res_W, (\eta\times\nu)\ind^W_Y} \cdot \inprod{(\eta\times\nu)\ind^W_Y\res_P, \sL}\\
		&= \sum_{\substack{\mu\in\cP(2^{k-1}),\\ r\in\{0,1\}}} a^\lambda_{\hook{2}{r},\mu} \cdot \inprod{\Chi(\mu;\phi_r)\res_P, \sL} + \sum_{\substack{\eta,\nu\in\cP(2^{k-1}),\\ \eta>\nu}} c^\lambda_{\eta,\nu} \cdot \sbc^{\eta}_{\cL(\hook{2^{k-1}}{z})} \cdot \sbc^{\nu}_{\cL(\hook{2^{k-1}}{z})}
	\end{align*}	
	\end{footnotesize}
	
	\noindent since $(\eta\times\nu)\ind^W_Y\res_P=(\eta\times\nu)\res^Y_B\ind^P$ by Mackey's Theorem \cite[\textsection 5]{Isaacs}. For $g_1,g_2\in Q=P_{2^{k-1}}$, we have by \cite[4.3.9]{JK} that $\Chi(\mu;\phi_r)(g_1,g_2;1) = \mu(g_1)\mu(g_2)$ and $\Chi(\mu;\phi_r)(g_1,g_2;(12)) = \mu(g_1g_2)\phi_r((12)) = \mu(g_1g_2)\cdot(-1)^r$, and we can similarly evaluate $\Chi\big(\cL(\hook{2^{k-1}}{z});\phi_i\big)$. Since $P=\{(g_1,g_2;\pi)\mid g_1,g_2\in Q,\ \pi\in P_2=\{1,(12)\}\}$ and $|P|=2|Q|^2$, we obtain
	\begin{small}
	\begin{align*}
		\inprod{\Chi(\mu;\phi_r)\res_P, \sL} &= \frac{1}{|P|}\sum_{g\in P} \Chi(\mu;\phi_r)(g) \cdot \overline{\Chi\big(\cL(\hook{2^{k-1}}{z});\phi_i\big)(g)}\\
		&= \frac{1}{2} \cdot \sbc^\mu_{\cL(\hook{2^{k-1}}{z})} \cdot \left( \sbc^\mu_{\cL(\hook{2^{k-1}}{z})} + (-1)^{r+i} \right).
	\end{align*}	
	\end{small}
	
	Next, note $a^\lambda_{(2),\mu} + a^\lambda_{(1^2),\mu}=c^\lambda_{\mu,\mu}$ by \eqref{eqn:a-c}, and $c^\lambda_{\mu,\mu}=1$ for $\mu\in\{(2^{k-1}),(2^{k-1}-1,1)\}$ and $c^\lambda_{\mu,\mu}=0$ for all other $\mu\in\cP(2^{k-1})$ by the Littlewood--Richardson rule. Moreover, $a^\lambda_{(2),(2^{k-1})}=0$ by \cite[Proposition 2.13]{LO} so $a^\lambda_{(1^2),(2^{k-1})}=1$. Also, $a^\lambda_{(1^2),(2^{k-1}-1,1)}\ge a^{(2^k-5,1)}_{(1^2),(2^{k-1}-2)}$ by \cite[Theorem 1.2]{dBPW21} and $a^{(2^k-5,1)}_{(1^2),(2^{k-1}-2)}=1$ by \cite[(5)]{W}, thus $a^\lambda_{(1^2),(2^{k-1}-1,1)}=1$ and $a^\lambda_{(2),(2^{k-1}-1,1)}=0$.
	
	By the Littlewood--Richardson rule, for $\eta>\nu\in\cP(2^{k-1})$ we have that $c^\lambda_{\eta,\nu}=1$ if $\eta=(2^{k-1})$ and $\nu=(2^{k-1}-j,j)$ for $j\in\{1,2,3\}$ or if $\eta=(2^{k-1}-1,1)$ and $\nu=(2^{k-1}2,2)$, and $c^\lambda_{\eta,\nu}=0$ otherwise. Hence
	\begin{equation}\label{eqn:ZXY}
		\sbc^{\lambda}_{\sL} = \sum_{\mu\in\{(2^{k-1}),(2^{k-1}-1,1)\}} \sX_\mu + \sum_{\nu\in\{(2^{k-1}-j,j)\mid j\in\{1,2,3\}\}} \sY_{(2^{k-1}),\nu} + \sY_{(2^{k-1}-1,1),(2^{k-1}-2,2)}
	\end{equation}
	where $\sX_\mu := \frac{1}{2} \cdot \sbc^\mu_{\cL(\hook{2^{k-1}}{z})} \cdot \left( \sbc^\mu_{\cL(\hook{2^{k-1}}{z})} + (-1)^{1+i} \right)$ and $\sY_{\eta,\nu} := \sbc^{\eta}_{\cL(\hook{2^{k-1}}{z})} \cdot \sbc^{\nu}_{\cL(\hook{2^{k-1}}{z})}$.
	Now, by \cref{thm: hook lins}, for $j\in\{0,1\}$ we have that $\sbc^{(2^{k-1}-j,j)}_{\cL(\hook{2^{k-1}}{z})}=1$ if $z=j$ and $\sbc^{(2^{k-1}-j,j)}_{\cL(\hook{2^{k-1}}{z})}=0$ otherwise. We also know the value of $\sbc^{(2^{k-1}-2,2)}_{\cL(\hook{2^{k-1}}{z})}$ from \cref{thm:almost-hook-sbcs}, and of $\sbc^{(2^{k-1}-3,3)}_{\cL(\hook{2^{k-1}}{z})}$ by the inductive hypothesis. In particular, the only non-zero values of $\sX_\mu$ and $\sY_{\eta,\nu}$ are:
	\begin{small}
	\begin{align*}
		\text{when }z=0:\quad& \sX_{(2^{k-1})} = \tfrac{1+(-1)^{1+i}}{2}, \quad \sY_{(2^{k-1}),(2^{k-1}-2,2)} = k-2, \quad \sY_{(2^{k-1}),(2^{k-1}-3,3)} = \tfrac{(k-1)(k-4)}{2};\\
		\text{when }z=1:\quad& \sX_{(2^{k-1}-1,1)} = \tfrac{1+(-1)^{1+i}}{2},\quad \sY_{(2^{k-1}-1,1),(2^{k-1}-2,2)} = k-3.
	\end{align*}	
	\end{small}
	
	\noindent The proof is then concluded by substituting these into \eqref{eqn:ZXY}, noting by \cref{rewrite lin} that $\Chi\big(\cL(\hook{2^{k-1}}{z});\phi_i\big) = \cL(\hook{2^k}{0})$, $\cL(\hook{2^k}{1})$, $\cL(\hook{2^k}{3})$, $\cL(\hook{2^k}{2})$ if $(z,i)=(0,0)$, $(0,1)$, $(1,0)$, $(1,1)$ respectively.
\end{proof}
}

\begin{theorem}\label{lins of 2^k}
	Let $k\in\N_{\geq 4}$ and $\lambda\in\Part(2^k)\setminus(\Hook(2^k)\cup\set{\ahook{2^k}{2^{k-1}-2}})$. Then $\abs{\Lin(\chi^\lambda\res_{P_{2^k}})} > 2$.
\end{theorem}
\begin{proof}
	If $\lambda$ is an almost hook which is not self-conjugate then $\abs{\Lin(\lambda\res_{P_{2^k}})}>2$ by \cref{cor: ah lins}. 
	Suppose $\lambda$ is not an almost hook. Recall from \cref{sec:rep-sn} that we may without loss of generality assume that $\lambda$ is wide or self-conjugate. We proceed with proof by induction on $k$. For $k=4$, the statement follows from direct computation. Now assume $k\ge 5$ and that the statement holds true for $k-1$. Let $\lambda\in\Part(2^k)\setminus(\Hook(2^{k})\cup\set{\ahook{2^k}{2^{k-1}-2}})$. 

	First suppose that $\lambda\not\in\cD_k$. Since $c^{\lambda'}_{\Delta(\lambda'),\Delta(\lambda')} = c^{\lambda}_{\Delta(\lambda')',\Delta(\lambda')'}$, it follows that $c^{\lambda}_{\mu,\mu}>0$ and $\abs{\Lin(\mu\res_{P_{2^{k-1}}})}>2$ for some $\mu\in\set{\Delta(\lambda),\Delta(\lambda')'}$ by \cref{split part}, \cref{annoying partitions} and the induction hypothesis. 
	Second, if $\lambda\in\cD_k\setminus\{(2^k-3,3),(2^3,1^{2^k-6})\}$, then there exist $\mu,\nu\in\Hook(2^{k-1})\cup\Ahook(2^{k-1})$ such that $c^\lambda_{\mu,\mu},c^\lambda_{\nu,\nu}>0$ and $\abs{\Lin(\mu\res_{P_{2^{k-1}}})\cup\Lin(\nu\res_{P_{2^{k-1}}})}>2$ by \cref{lem:fixing spc cases}. In both cases, it follows that $\abs{\Lin(\lambda\res_{P_{2^k}})}>2$ by \eqref{eq:find lins}.
	
	{\blue Finally, suppose that $\lambda\in\{(2^k-3,3),(2^3,1^{2^k-6})\}$. Since $\abs{\Lin(\chi^\lambda\res_{P_n})} = \abs{\Lin(\chi^{\lambda'}\res_{P_n})}$, it suffices to show that $|\Lin(\chi^{(2^k-3,3)}\down_{P_{2^k}})|>3$. This follows immediately from \cref{lem:2k-3}.}
\end{proof}

  We are now ready to prove \cref{thm:main-p=2}, which determines precisely those $\chi\in\Irr(S_n)$ for which $\abs{\Lin(\chi\res_{P_n})}=2$ where $P_n\in\Syl_2(S_n)$.

  \begin{proof}[Proof of \cref{thm:main-p=2}]
		Let $n\in\N$, $P_n\in\Syl_2(S_n)$ and $\lambda\in\Part(n)$. Recall from \cref{sec:rep-sn} that we may assume $\lambda$ is wide or self-conjugate when counting $\abs{\Lin(\chi\res_{P_n})}$.
		It follows from \cref{cor: ah lins} that the partitions $\lambda$ described in case (1) satisfy $\abs{\Lin(\lambda\res_{P_n})}=2$ and $\inprod{\lambda\res_{P_n},\psi}=1$ for all $\psi\in\Lin(\lambda\res_{P_n})$. Case (2) follows from \cref{thm: hook lins} and the Littlewood--Richardson rule, and case (3) follows from direct computation.
		
		Now assume that $\lambda$ is not any of the partitions described in (1), (2) or (3). We proceed to show that $\abs{\Lin(\lambda\res_{P_n})}\neq 2$. Firstly note for all $n\in\N$ that if $\lambda\in\set{(n),(1^n)}$, then $\abs{\Lin(\lambda\res_{P_n})}=1$. 
		Let $m,k\in\N_0$ be such that $m<2^k$ and $n=2^k+m$. If $k<4$, then the statement follows from direct computation. If $m=0$, then the statement follows from \cref{thm: hook lins} and \cref{lins of 2^k}. Left to consider is when $k\ge 4$ and $m\geq 1$. We will consider the two cases $\lambda\in\Hook(n)$ and $\lambda\in\Part(n)\setminus\Hook(n)$ separately.
		Firstly, suppose that $\lambda\in\Part(n)\setminus\Hook(n)$. Then there exists some $\mu\in\Part(2^k)\setminus(\Hook(2^k)\cup\{\ahook{2^{k}}{2^{k-1}-2}\})$ such that $\mu\subseteq\lambda$. Since $\abs{\Lin(\mu\res_{P_{2^k}})}>2$ by \cref{lins of 2^k} it follows that $\abs{\Lin(\lambda\res_{P})}>2$ by \cref{lem: lins lower bound}.	
		Now, suppose that $\lambda\in\Hook(n)$, so $\lambda=\hook{n}{x}$ for some $x\in\set{1,\dots,\lfloor\frac{n-1}{2}\rfloor}$. Since $\lambda$ is none of the partitions in (1)--(3), we have that $m\geq 2$. 
		
		If $\lambda\ne\hook{n}{1}$, then since $m\ge 2$ there exist by the Littlewood--Richardson rule distinct $a_1,a_2,a_3\in\set{0,1,\dotsc,2^{k-1}-1}$ such that $\hook{2^k}{a_i}\subseteq\lambda$ for all $i\in\set{1,2,3}$. Namely, we may take $a_1=\max\set{0,2^k-\lambda_1}$, $a_2=a_1+1$ and $a_3=a_1+2$.
		So there exist $\nu_{(1)},\nu_{(2)},\nu_{(3)}\in\Part(m)$ such that $\lambda\res_{S_{2^k}\times S_m}$ contains $\hook{2^k}{a_i}\times\nu_{(i)}$ as an irreducible constituent for all $i\in\set{1,2,3}$. Since $P_{2^k}\times P_m\cong P_n$ by \eqref{eqn:Pn}, we deduce that there exist $\theta_{(1)},\theta_{(2)},\theta_{(3)}\in\Lin(P_m)$ such that $\lin(\hook{2^k}{a_i})\times\theta_{(i)}\in\Lin(\lambda\res_{P_n})$ for all $i\in\set{1,2,3}$ by \cref{prop: at least one lin} and \cref{thm: hook lins}. But the $a_i$ are distinct, so $\abs{\Lin(\lambda\res_{P_n})}\ge 3$.
		
		Finally, suppose $\lambda=\hook{n}{1}$. Write $m=2^t+u$ with $0\le u<2^t$, so $t\ge 1$ since $m\ge 2$. That is, $n=2^k+2^t+u$ and $k>t\ge 1$. Then the Littlewood--Richardson rule gives that 
		\[ \hook{2^k}{1}\times\hook{2^t}{0}\times\hook{u}{0},\quad \hook{2^k}{0}\times\hook{2^t}{0}\times\hook{u}{0},\quad \hook{2^k}{0}\times\hook{2^t}{1}\times\hook{u}{0} \]
		are irreducible constituents of $\lambda\res_{S_{2^k}\times S_{2^t}\times S_u}$. Since $P_{2^k}\times P_{2^t}\times P_u\cong P_n$ by \eqref{eqn:Pn}, we deduce that $\lin(\hook{2^k}{1})\times\triv_{P_{2^t}}\times\triv_{P_u}$, $\triv_{P_n}$ and $\triv_{P_{2^k}}\times\lin(\hook{2^t}{1})\times\triv_{P_u}$ are three distinct elements of $\Lin(\lambda\res_{P_n})$ by \cref{thm: hook lins}, noting that $\lin(\hook{2^k}{0})=\triv_{P_{2^k}}$ and similarly for $t$.
  \end{proof}

\section{Odd primes}\label{sec:odd}

Let $n\in\N_0$ and $p$ be a prime. In \cref{thm:main-p=2} we classified the set of $\chi\in\Irr(S_n)$ for which $\abs{\Lin(\chi\res_{P_n})}=p$ for $p=2$, where $P_n\in\Syl_p(S_n)$. We now proceed to complete the classification for all odd primes. Many of the important results used in the previous sections cannot be extended to the case of odd $p$, but instead we can use certain results on the positivity of Sylow branching coefficients in this case from \cite{GL2,L}. 

\begin{definition}
	Let $n,t\in\N$ and $\psi\in\Irr(P_n)$.
	\begin{itemize}
		\item[i)] We define $\Omega(\psi):=\set{\chi\in\Irr(S_n)\mid \sbc^\chi_\psi>0}$.
		\item[ii)] Let $\cB_n(t):= \set{\lambda\in\Part(n) \mid \lambda_1, \ell(\lambda) \leq t}$.
		\item[iii)] Let $m(\psi):=\max\set{t\in\N \mid \cB_n(t) \subseteq \Omega(\psi)}$ and $M(\psi):=\min\set{t\in\N \mid \Omega(\psi)\subseteq \cB_n(t)}$.  
	\end{itemize} 	
\end{definition}

Hence we have for each $\psi\in\Lin(P_n)$ that $\cB_n(m(\psi))\subseteq \Omega(\psi)\subseteq \cB_n(M(\psi))$. On the one hand, if $\lambda\in\cB_n(m(\psi))$ then $\psi\in\Lin(\chi^\lambda\res_{P_n})$. On the other hand, if $\lambda\not\in\cB_{n}(M(\psi))$ then $\psi\not\in\Lin(\chi^\lambda\res_{P_n})$. Hence, knowledge of the values of $m(\psi)$ and $M(\psi)$ can be used to bound the size of $\Lin(\chi^\lambda\res_{P_n})$ for $\lambda\in\Part(n)$.

We first record what happens in the case $n=p$. Recall that $P_p\in\Syl_p(S_p)$ is cyclic and that we denote $\Irr(P_p):=\set{\phi_0,\phi_1,\dots,\phi_{p-1}}$, where $\phi_0$ is the trivial character.
\begin{lemma}\cite[Lemma 4.2]{GL2}\label{Omega k=1}
	Let $p$ be an odd prime, let $P\in\Syl_p(S_p)$ and $\lambda\in\Part(p)$. Then $\phi_0$ is a constituent of $\chi^\lambda\res_P$ if and only if $\lambda\notin\{(p-1,1),(2,1^{p-2})\}$. For each $i\in\{1,2,\dotsc,p-1\}$, $\phi_i$ is a constituent of $\chi^\lambda\res_P$ if and only if $\lambda\in\cB_p(p-1)$.
\end{lemma}

Before we proceed, we need to introduce some more notation. 
Let $i,j\in\set{0,1,\dots,p-1}$ and $k\in\N$. We define 
\begin{equation}\label{eqn:linears}
	\psi_i:=\Chi(\triv_{P_{p^{k-1}}};\phi_i),\quad \text{and for } k\geq 2,\ \text{let }\psi_{i,j}:=\Chi(\Chi(\triv_{P_{p^{k-2}}};\phi_i);\phi_j).
\end{equation}
\noindent In particular, $\psi_{0,j}=\psi_{j}$, $\psi_{0}=\triv_{P_{p^k}}$ and if $k=1$ then $\psi_i =\phi_i$. The next lemma follows from \cite[Lemma 4.3, Theorem 4.5, Theorem 4.9]{GL2} and \cite[Theorem 2.18]{L}.

\begin{lemma}\label{compute Mm}
	Let $p$ be an odd prime and $k\in\N_{\geq 2}$. 
		\begin{enumerate}[label=\roman*)]
			\item Suppose that $\psi\in\Lin(P_{p^k})\setminus\set{\psi_0,\dots,\psi_{p-1}}$. Then $M(\psi)\leq p^k-p$, with equality if and only if $\psi=\psi_{i,j}$ for some $i\in\set{1,2,\dots,p-1}$, $j\in\set{0,1,\dots,p-1}$.
			\item Suppose that $(p,k)\neq (3,2)$ and $i\in\set{1,2,\dots,p-1}$. Then $m(\psi_i)=M(\psi_i)=p^k-1$ and $m(\psi_{i,0})=p^k-(p+1)$.
		\end{enumerate}
\end{lemma}

We can now explicitly compute $\Lin(\chi^\lambda\res_{P_{p^k}})$ for certain $\lambda\in\Part(p^k)$ by using these values of $m(\psi)$ and $M(\psi)$ where $\psi\in\Lin(P_{p^k})$.

\begin{lemma}\label{lem: small lins}
	Let $k\in\N$ and $p$ be an odd prime. Then $\Lin(\chi^{(p^k-1,1)}\res_{P_{p^k}}) = \set{\psi_1,\psi_2,\dots,\psi_{p-1}}$, and $\Lin(\chi^\lambda\res_{P_{p^k}})=\set{\psi_0,\psi_1,\dots,\psi_{p-1}}$ for all $\lambda\in\cB_{p^k}(p^k-2)\setminus\cB_{p^k}(p^k-p)$.
\end{lemma}

\begin{proof}
	For $k=1$, the statement follows from \cref{Omega k=1}. Now suppose $k\geq 2$ and let $\lambda\in\cB_{p^k}(p^k-1)\setminus\cB_{p^k}(p^k-p)$. Let $i\in\set{1,2,\dots,p-1}$. Then $m(\psi_i) = p^k-1$ by \cref{compute Mm} ii). Since $\lambda \in\cB_{p^k}(p^k-1)$, it follows that $\psi_i\in\Lin(\lambda\res_{P_{p^k}})$. Furthermore, $\psi_0 \in\Lin(\lambda\res_{P_{p^k}})$ if and only if $\hook{p^k}{1}\notin\{\lambda,\lambda'\}$ by \cite[Theorem A]{GL1}. Next, let $\psi\in\Lin(P_{p^k})\setminus\set{\psi_0,\psi_1,\dots,\psi_{p-1}}$. Then $M(\psi)\leq p^k-p$  by \cref{compute Mm} i), and since $\lambda\not\in\cB_{p^k}(p^k-p)$ it follows that $\psi\not\in\Lin(\lambda\res_{P_{p^k}})$. 
\end{proof}

\begin{remark}
	In fact, when $k\in\N$ and $p$ is an odd prime, then $\sbc^{(p^k-1,1)}_{\psi_i}=1$ 
	for all $i\in\{1,\dotsc,p-1\}$. This is clear for $k=1$. To see this for $k\ge 2$, let $Q:=P_{p^{k-1}}$, $B=Q^{\times p}$ (so that $P:=P_{p^k}=Q\wr P_p = B\rtimes P_p$) and $Y:=S_{p^{k-1}}^{\times p}$. The Littlewood--Richardson rule gives
	\[ \inprod{\chi^{(p^k-1,1)}\res_B,\triv_B} = \inprod{\chi^{(p^k-1,1)}\res_Y\res_B,\triv_B} = p-1 \]
	since $\chi^{(p^k-1,1)}\res_Y = (p-1)\cdot\triv_Y + \sum (\phi\times\triv\times\cdots\times\triv)$ where the sum runs over all $p$ cyclic permutations of the factors and $\phi:=\chi^{(p^{k-1}-1,1)}$ and $\triv:=\triv_{S_{p^{k-1}}}$, noting that $\sbc^{(p^{k-1}-1,1)}_{\triv_Q}=0$ by \cite[Theorem A]{GL1}. But also 
	\[ \inprod{\chi^{(p^k-1,1)}\res_B,\triv_B} \ge \sum_{i=1}^{p-1} \sbc^{(p^k-1,1)}_{\psi_i}\cdot \inprod{\psi_i\res_B,\triv_B} = \sum_{i=1}^{p-1} \sbc^{(p^k-1,1)}_{\psi_i} \ge p-1 \]
	since $\psi_i\res_B = \triv_B$ and $\sbc^{(p^k-1,1)}_{\psi_i}>0$ from \cref{lem: small lins}. Thus $\sbc^{(p^k-1,1)}_{\psi_i}=1$ for all $i\in\{1,\dotsc,p-1\}$. 
	
	We remark that a similar argument can be used to show that $\chi^{(p^k-1,1)}\res_P$ is multiplicity-free. In other words, if $\phi\in\Irr(P)$ satisfies $\sbc^{(p^k-1,1)}_{\phi}>0$ then $\sbc^{(p^k-1,1)}_{\phi}=1$. This follows, using induction on $k$, from the decomposition of $\chi^{(p^k-1,1)}\res_Y$ and $\sbc^{(p^{k-1}-1,1)}_{\triv_Q}=0$, and the fact that $\Irr(P\mid \tau\times\triv_Q^{\times(p-1)})=\{\psi_\tau:=(\tau\times\triv_Q^{\times(p-1)})\ind^P_B\}$ whenever $\tau\in\Irr(Q)\setminus\{\triv_Q\}$, which together imply that $\sbc^{(p^k-1,1)}_{\psi_\tau}=\sbc^{(p^{k-1}-1,1)}_{\tau}$.
\end{remark}

\begin{proof}[Proof of \cref{thm:main-odd-p}]
	Recall from \cref{sec:rep-sn} that we may without loss of generality assume that $\lambda$ is wide or self-conjugate. 
	
	\noindent (1) Suppose that $n=p^k$ for some $k\in \N$. If $(p,k)=(3,2)$ then the statement follows from \cite[p.~20]{L}, so now suppose $(p,k)\ne (3,2)$.
	If $\lambda\not\in\cB_n(n-p)$, then the statement follows from \cref{lem: small lins}. 
	
	Now suppose that $\lambda\in\cB_n(n-p)$. If $k=1$ then $\cB_n(n-p)=\varnothing$ and so we may assume that $k\geq 2$. 	
	If $j\neq 0$, then $m(\psi_j)=n-1$ by \cref{compute Mm} ii). Furthermore, $\psi_0\in\Lin(\lambda\res_{P_n})$ by \cite[Theorem A]{GL1}. Hence $\set{\psi_0,\psi_1,\dots,\psi_{p-1}}\subseteq\Lin(\lambda\res_{P_n})$. 
	We now show there exists some $\varnothing\ne L\subseteq\Lin(P_n)\setminus\set{\psi_0,\psi_1,\dots,\psi_{p-1}}$ such that $L\subseteq\Lin(\lambda\res_{P_n})$. 
	\begin{itemize}
		\item First, suppose that $\lambda\in\cB_n(n-(p+1))$. Then $m(\psi_{i,0})=n-p$ for $i\neq 0$ by \cref{compute Mm} ii), and so $\set{\psi_{1,0},\dots,\psi_{p-1,0}}\subseteq \Lin(\lambda\res_{P_n})$. 
		
		\item Next, suppose that $\lambda\in\set{(n-p)\sqcup \mu \mid \mu\in\Part(p)}$. If $i\neq 0$ then $M(\psi_{i,j})=n-p$ by \cref{compute Mm} i). If $\mu\in\set{(p),(1^p)}$ then $\set{\psi_{1,0},\dots,\psi_{p-1,0}}\subseteq\Lin(\lambda\res_{P_n})$ by \cite[Lemma 4.6]{GL2} and \cref{Omega k=1}. 
		
		\item Finally, if $\mu\not\in\set{(p),(1^p)}$, then $\set{\psi_{i,j}\mid i,j\in\set{1,\dots,p-1}}\subseteq\Lin(\lambda\res_{P_n})$ by \cite[Lemma 4.6]{GL2} and \cref{Omega k=1}.
	\end{itemize} 
	Thus $\abs{\Lin(\lambda\res_{P_n})}>p$ as desired.
	
	\smallskip

	\noindent (2) Now suppose that $n\neq p^k$ for any $k\in\N$. We split (2) into cases (I) and (II) as follows.
		
	\noindent \textbf{Case (I): $n=p^k+i$ for some $i\in\set{1,2,\dots,p-1}$.} 
	First, suppose $\lambda\in\cB_n(n-1)\setminus\cB_n(n-p)$. 
	By the Littlewood--Richardson rule, we have that
	\[ \lambda\res_{P_n} = \sum_{\substack{\mu\in\Part(p^k)\setminus\cB_{p^k}(p^k-p),\\ \nu\in\Part(i)}} c^\lambda_{\mu,\nu} \cdot \mu\res_{P_{p^k}}\times \nu\res_{P_1^{\times i}} = \sum_{\substack{\mu\in\Part(p^k)\setminus\cB_{p^k}(p^k-p),\\ \nu\in\Part(i)}} c^\lambda_{\mu,\nu}\cdot\chi^\nu(1)\cdot \mu\res_{P_{p^k}}\times\triv_{P_1}^{\times i}. \]
	By considering $\Lin(\mu\res_{P_{p^k}})$ for such $\mu$, it follows from \cref{lem: small lins} that $\Lin(\lambda\res_{P_n})\subseteq \set{\psi_j\times\triv_{P_1}^{\times i} \mid j\in\set{0,\dots,p-1}}$. We show that these two sets are in fact equal {\blue whenever $(p,n,\lambda)\ne\big(3,4,(2^2)\big)$}.
	\begin{itemize}
		\item If $\lambda=(n-1,1)$, then $(p^k-1,1)\subseteq\lambda$ and so $\set{\psi_j\times \triv_{P_1}^{\times i} \mid j\in\set{1,\dots,p-1}}\subseteq \Lin(\lambda\res_{P_n})$ by \cref{lem: small lins}. Moreover, $(p^k)\subseteq\lambda$ so $\psi_0\times\triv_{P_1}^{\times i}\in\Lin(\lambda\res_{P_n})$ also.
		\item If $\lambda\ne(n-1,1)$: {\blue if $(p,k)\ne(3,1)$} then there exists $\nu\in\cB_{p^k}(p^k-2)\setminus\cB_{p^k}(p^k-p)$ such that $\nu\subseteq\lambda$.		
		Hence $\set{\psi_j\times\triv_{P_1}^{\times i} \mid j\in\set{0,\dots,p-1}} \subseteq \Lin(\lambda\res_{P_n})$ by \cref{lem: small lins}.
		{\blue If $(p,k)=(3,1)$ then $n\in\{4,5\}$. When $n=4$ then $\Lin((2^2)\down_{P_4})=\Irr(P_4)\setminus\{\triv_{P_4}\}$. When $n=5$ then $\Lin(\lambda\down_{P_5})=\Irr(P_5)$ for all $\lambda\in\{(3,2),(3,1^2),(2^2,1)\}$.}
	\end{itemize}
	Thus $\abs{\Lin(\lambda\res_{P_n})}=p$ for all $\lambda\in\cB_n(n-1)\setminus\cB_n(n-p)$, {\blue unless $(p,n,\lambda)=\big(3,4,(2^2)\big)$ in which case $\abs{\Lin(\lambda\res_{P_n})}=p-1$.}
	
	Next, suppose that $\lambda\in\cB_n(n-p)$. Then $\lambda=(n-p-t)\sqcup \mu$ for some $t\in\N_0$ and $\mu\in\Part(p+t)$. First suppose $k\ge 2$. Since $\abs{\mu}\geq p$, there exists some $\nu\in\cB_{p^k}(p^k-p)$ such that $\nu\subseteq\lambda$, and $\nu\neq(3^3)$ if $(p,k)=(3,2)$. It follows from \cref{thm:main-odd-p} (1) and \cref{lem: lins lower bound} that $\abs{\Lin(\lambda\res_{P_n})}>p$. 
	{\blue Otherwise, $k=1$, so $\lambda\in\cB_n(i)$. Let $\Psi_j:=\psi_j\times\triv_{P_1}^{\times i}$ for $j\in\{0,1,\dotsc,p-1\}$. For $j\in\{1,\dotsc,p-1\}$ we have $i\le n-1=m(\Psi_j)$ by \cite[Theorem 2.11]{GL2} and \cite[Table 2]{L}. Hence $\cB_n(i)\subseteq \cB_n(m(\Psi_j))$ and so $\set{\Psi_1,\dots,\Psi_{p-1}}\subseteq \Lin(\lambda\res_{P_n})$. 
	Furthermore, $\Psi_0\in \Lin(\lambda\res_{P_n})$ by \cite[Theorem A]{GL1}. Hence $\Lin(\lambda\res_{P_n}) = \Lin(P_n)$ and so $\abs{\Lin(\lambda\res_{P_n})}=p$. }
	
	\smallskip
	
	\noindent \textbf{Case (II): $n\in\set{p^k+p,\dots,p^{k+1}-1}$.}
	Suppose $n=\sum_{i\ge 0}b_ip^i$ in $p$-adic expansion. Let $r=\sum_i b_i$ and let $a_1\le\dotsc\le a_r\in\N_0$ be such that $\#\set{j\in\set{1,2,\dotsc,r} \text{ s.t. } a_j=i}=b_i$ for each $i\ge 0$, so $n=p^{a_1}+\cdots+p^{a_r}$. Let $t=\min\{i\mid a_i>0\}$, so $t<r$ since $n\ge p^k+p$. 
	
	Let $X=S_{p^{a_1}}\times\cdots\times S_{p^{a_r}}$ and observe that $P_n\cong P_{p^{a_1}}\times\cdots\times P_{p^{a_r}}\le X$. For $l\in\set{1,2,\dotsc,r}$ and $j\in\set{0,1,\dotsc,p-1}$, let 
	\[ \psi_j^{(l)}:=(\triv,\dotsc,\psi_j,\dotsc,\triv) \in \Lin(P_{p^{a_1}})\times\cdots\times\Lin(P_{p^{a_r}}) = \Lin(P_n) \]
	where $\psi_j\in\Lin(P_{p^{a_l}})$ is defined as in \eqref{eqn:linears}.
	\begin{itemize}
		\item If $\lambda=(n-1,1)$, then $\triv\times\cdots\times \hook{p^{a_t}}{1}\times\cdots\times\triv$ and $\triv\times\cdots\times\triv\times\hook{p^{a_r}}{1}$ are distinct irreducible constituents of $\lambda\res_X$ by the Littlewood--Richardson rule. Hence $\psi_j^{(t)},\psi_j^{(r)}\in\Lin(\lambda\res_{P_n})$ for all $j\in\set{1,2,\dotsc,p-1}$ by \cref{lem: small lins}, and so $\abs{\Lin(\lambda\res_{P_n})}\ge 2p-2>p$.
		
		\item If $\lambda=(\lambda_1,\mu)$ where $2\le|\mu|<p$ then $\triv\times\cdots\times\triv\times\hook{p^{a_t}}{1}\times\triv\times\cdots\times\triv\times\hook{p^{a_r}}{1}$ is an irreducible constituent of $\lambda\res_X$. Hence $\psi_j^{(t)}\cdot\psi_i^{(r)}\in\Lin(\lambda\res_{P_n})$ for all $i,j\in\set{1,2,\dotsc,p-1}$ by \cref{lem: small lins}, and so $\abs{\Lin(\lambda\res_{P_n})}\ge (p-1)^2>p$.
		
		\item Otherwise, $\lambda=(\lambda_1,\mu)$ where $|\mu|\ge p$. If $k=1$, then $n=ap+b$ where $2\le a<p$ and $0\le b<p$. By \cite[Proposition 3.2]{GL2}, $\cB_p(p-1)^{\ast a}\ast\cB_1(1)^{\ast b} = \cB_n(n-a) \supseteq \cB_n(n-p)\ni\lambda$, so $\lambda\res_X$ has an irreducible constituent of the form $\gamma_1\times\cdots\times\gamma_a\times(1)^{\times b}$ where $\gamma_i\in\cB_p(p-1)$ for all $i$. But $\Lin(\gamma_i\res_{P_p})=\set{\phi_1,\dotsc,\phi_{p-1}}$ by \cref{Omega k=1}, so $\phi_{i_1}\times\cdots\times\phi_{i_a}\times\triv\times\cdots\times\triv\in\Lin(\lambda\res_{P_n})$ for all $i_j\in\set{1,\dotsc,p-1}$. Hence $\abs{\Lin(\lambda\res_{P_n})}\ge(p-1)^a>p$.
		
		If $k\ge 2$ then there exists $\nu\in\cB_{p^k}(p^k-p)$ such that $\nu\subseteq\lambda$ (and $\nu\ne(3^3)$ if $(p,k)=(3,2)$), whence $\abs{\Lin(\lambda\res_{P_n})}\ge\abs{\Lin(\nu\res_{P_{p^k}})}>p$ by \cref{lem: lins lower bound} and \cref{thm:main-odd-p} (1).
	\end{itemize}
	Thus $\abs{\Lin(\lambda\res_{P_n})}>p$ for all $\lambda$ considered in Case (II), as desired.
\end{proof}

{\blue
\begin{remark}\label{rem:identify-few-linears-odd}
	Let $p$ be an odd prime.
	When $p=3$ then $\Lin((3^3)\down_{P_9})=\{\Chi(\phi_i;\phi_0)\mid i\in\{0,1,2\} \}$ by direct computation. Together with \cref{lem: small lins}, this explicitly identifies all of the elements of $L:=\Lin(\lambda\down_{P_n})$ when $|L|\le p$, when $n=p^k$ for some $k\in\N$ (i.e.~those in case (1) of \cref{thm:main-odd-p}).
	
	Next, we identify the elements of $L$ when $|L|\le p$ and $n$ is not a power of $p$, i.e.~ those in case (2) of \cref{thm:main-odd-p}. 
	We note that $\Lin((2^2)\down_{P_4})=\Irr(P_4)\setminus\{\triv_{P_4}\}$ when $p=3$. 
	Now let $n=p^k+i$ for some $k\in\N$ and $i\in\{1,2,\ldots,p-1\}$.
	If $k=1$ and $\lambda\in\cB_n(n-p)$ then $|L|=p$ implies $L=\Irr(P_n)$ as $P_n\cong C_p$. Otherwise, $k\in\N$ and $\lambda\in\cB_n(n-1)\setminus\cB_n(n-p)$ and $(p,n,\lambda)\ne(3,4,(2^2))$, in which case $L=\{\Chi(\triv_{P_{p^{k-1}}};\phi_j) \times\triv_{P_i} \mid 0\le j\le p-1 \}$ from the proof of \cref{thm:main-odd-p} (2) Case (I).
\end{remark} }

In this paper, we have thus far focused on irreducible characters of symmetric groups whose restrictions to a Sylow subgroup $P$ possess small numbers of distinct linear constituents. At the other end of the spectrum, there are in general many partitions $\lambda$ for which $|\Lin(\chi^\lambda\down_P)|$ takes the maximum possible value. We give an example to illustrate this observation.

\begin{example}
	Let $k\in\N$. Then $P\in\Syl_p(S_{p^k})$ has exactly $p^k$ linear characters since $P/P'\cong(C_p)^k$. From \cite[Theorem 2.11]{GL2}, any $\lambda\in\cB_{p^k}(p^k-p^{k-1}-p^{k-2})$ has $|\Lin(\chi^\lambda\down_P)|=p^k$ when $p\ge 5$. Similarly, $\lambda\in\cB_{3^k}(\tfrac{3^k}{2})$ implies $|\Lin(\chi^\lambda\down_P)|=3^k$ when $k>4$ by \cite[Theorem 2.18]{L}.
	It would be interesting to investigate the possible range of values of $|\Lin(\chi^\lambda\down_P)|$ when $p=2$.
\end{example}

\subsection*{Acknowledgements}
We thank the referee for their helpful comments and suggestions. We are grateful to Gabriel Navarro for his helpful comments on a previous version of this article. The first author thanks Eugenio Giannelli for helpful conversations on this topic and the Algebra research group at Universit\`a degli Studi di Firenze for their hospitality during which part of this work was done. We also thank David Craven for providing useful data for this project, and the School of Mathematics at the University of Birmingham and EPSRC for funding this research.


\end{document}